\documentclass[leqno,12pt]{article}
\usepackage{amsmath,amssymb}
\usepackage{color}
\usepackage{epsfig}

\newtheorem{Theorem}{\hspace{\parindent}\bf Theorem}[section]
\newtheorem{Lemma}{\hspace{\parindent}\bf Lemma}[section]
\newtheorem{Proposition}{\hspace{\parindent}\bf Proposition}[section]
\newtheorem{Corollary}{\hspace{\parindent}\bf Corollary}[section]

\newtheorem{Definition}{\hspace{\parindent}\bf Definition}[section]

\newcommand{\qed}{\hfill$\square$\vspace{0.3cm}}

\newcommand{\Tr}{\mathop{\rm Tr}}
\newcommand{\E}{\mathop{\rm E}}

\begin{document}

\title{\textbf{A general fractional \\ porous medium equation }}
\author{ by \\
Arturo de Pablo, Fernando Quir\'{o}s, \\ Ana Rodr\'{\i}guez, and Juan Luis V\'{a}zquez}

\maketitle

\begin{abstract}
We develop a theory of existence and uniqueness for the following
porous medium equation  with fractional diffusion,
$$
\left\{
\begin{array}{ll}
\dfrac{\partial u}{\partial t} + (-\Delta)^{\sigma/2}
(|u|^{m-1}u)=0, & \qquad  x\in\mathbb{R}^N,\; t>0,
\\ [8pt]
u(x,0) = f(x), & \qquad x\in\mathbb{R}^N.%
\end{array}
\right.
$$
We consider data $f\in L^1(\mathbb{R}^N)$ and all exponents $0<\sigma<2$ and $m>0$.
Existence and uniqueness of a weak solution is established for $m>
m_*=(N-\sigma)_+ /N$,  giving rise to an $L^1$-contraction semigroup. In addition, we
obtain the main qualitative properties of these solutions. In the lower range
$0<m\le m_*$ existence and uniqueness of solutions with good properties happen
under some restrictions, and the properties are different from the case above
$m_*$.  We also study the dependence of solutions on $f,m$ and
$\sigma$. Moreover, we consider the above questions for the problem posed in a
bounded domain.
\end{abstract}

\vskip 1cm

%
\noindent{\makebox[1in]\hrulefill}\newline

\newpage

\begin{center}\begin{minipage}{12cm}{\tableofcontents}\end{minipage}\end{center}

\

\newpage

\section{Introduction}\label{sect-intro}
\setcounter{equation}{0}

The aim of this paper is to develop a theory of existence and
uniqueness, as well as to obtain the main qualitative properties,
for a family of nonlinear fractional diffusion equations of porous
medium type. More specifically, we consider the Cauchy problem
\begin{equation}  \label{eq:main}
\left\{
\begin{array}{ll}
\dfrac{\partial u}{\partial t} + (-\Delta)^{\sigma/2}
(|u|^{m-1}u)=0, & \qquad  x\in\mathbb{R}^N,\; t>0,
\\ [4mm]
u(x,0) = f(x), & \qquad x\in\mathbb{R}^N.%
\end{array}
\right.
\end{equation}
We take initial data $f\in L^1(\mathbb{R}^N)$, which is a standard assumption
in diffusion problems,  with no sign restriction.  As for the exponents, we
consider the fractional exponent range $0<\sigma<2$, and take porous medium
exponent $m>0$. In the limit $\sigma\to 2$ we recover the standard
Porous Medium Equation (PME)
$$
\dfrac{\partial u}{\partial t}-\Delta (|u|^{m-1}u)=0,
$$
which is a basic model for nonlinear and degenerate  diffusion,
having now a well-established theory \cite{Vazquez}.

The nonlocal operator $(-\Delta)^{\sigma/2}$, known as the Laplacian
of order $\sigma$, is defined for any function $g$ in the  Schwartz
class through the Fourier transform: if $(-\Delta)^{\sigma/2}  g=h$,
then
\begin{equation}
\widehat{h}\,(\xi)=|\xi|^\sigma \,\widehat{g}(\xi).
\label{def-fourier}
\end{equation}
If $0<\sigma<2$, we can also use the representation by means of an
hypersingular kernel,
\begin{equation}
(-\Delta)^{\sigma/2}  g(x)= C_{N,\sigma }\mbox{
P.V.}\int_{\mathbb{R}^N} \frac{g(x)-g(z)}{|x-z|^{N+{\sigma} }}\,dz,
\label{def-riesz}
\end{equation}
where $C_{N,\sigma
}=\frac{2^{\sigma-1}\sigma\Gamma((N+\sigma)/2)}{\pi^{N/2}\Gamma(1-\sigma/2)}$
is a normalization constant, see for example \cite{Landkof}.
There is another classical way of defining the fractional powers of
a linear self-adjoint nonnegative operator, in terms of the
associated semigroup, which in our case reads
\begin{equation}
\displaystyle(-\Delta)^{\sigma/2}
g(x)=\frac1{\Gamma(-\frac{\sigma}2)}\int_0^\infty\left(e^{t\Delta}g(x)-g(x)\right)\frac{dt}{t^{1+\frac{\sigma}2}}.
\label{laplace}\end{equation} It is easy to check that the symbol of
this operator is again $|\xi|^\sigma$. The advantage of this approach is that it gives a natural way of defining the problem in a bounded domain, by means of the spectral characterization of
the semigroup $e^{t\Delta}$, see \cite{Stinga-Torrea}. The problem posed in a bounded domain is also studied in this paper.

Equations of the form \eqref{eq:main} can be seen as
fractional-diffusion versions of the PME.  Though our paper is aimed
at providing a sound mathematical theory for this evolution equation
and the nonlinear semigroups generated by Problem \eqref{eq:main}
and the problem posed on bounded domains, we mention that  the
equation appears as a model in statistical mechanics \cite{Jara},
and the linear counterpart in \cite{JKOlla}. We also want to point
out that there are other natural options of nonlinear, possibly
degenerate fractional-diffusion evolutions under current
investigation. Thus, the papers \cite{Caffarelli-Vazquez09},
\cite{Caffarelli-Vazquez10} consider the following fractional
diffusion PME  \  $\partial_t u - \nabla \cdot (u\nabla
(-\Delta)^{-s/2} u)=0$. It has very different properties from the
ones we derive for \eqref{eq:main}.  The standard PME (with $m=2$)
is recovered in such model for $s=0$.   A more detailed discussion
on these issues is contained in the survey paper \cite{Vazquez11}.

These two kinds of equations can be also viewed as nonlinear versions of the
linear fractional diffusion equation obtained for $m=1$, which has
the integral representation
\begin{equation}
u(x,t)=\int_{\mathbb{R}^N}K_\sigma(x-z,t)f(z)\,dz, \label{lineal}
\end{equation}
where $K_\sigma$ has Fourier transform  $\widehat
K_\sigma(\xi,t)=e^{-|\xi|^\sigma t}$. This means that, for
$0<\sigma<2$, the kernel $K_\sigma$ has the form
$K_\sigma(x,t)=t^{-N/\sigma}F(|x|\,t^{-1/\sigma})$ for some profile
function $F$ that is positive and decreasing, and behaves at
infinity like $F(r)\sim r^{-(N+\sigma)}$ \cite{Blumenthal-Getoor}.
When $\sigma=1$, $F$ is explicit; if $\sigma=2$ the function $K_2$
is the Gaussian heat kernel. The linear model has been well studied
by probabilists, since the fractional Laplacians of order $\sigma$
are infinitesimal generators of stable L\'{e}vy processes
\cite{Applebaum}, \cite{Bertoin}. However, an integral
representation of the evolution like \eqref{lineal} is not available
in the nonlinear case, thus motivating our work.

In a previous article \cite{dePablo-Quiros-Rodriguez-Vazquez} we
studied  Problem \eqref{eq:main}   for the particular case
$\sigma=1$. The key tool there was the well-known representation of
the half-Laplacian in terms of the Dirichlet-Neumann operator, which
allowed us to formulate the nonlocal problem in terms of a local one
(i.\,e., involving only derivatives and not integral operators). For
$\sigma\ne 1$, Caffarelli and Silvestre \cite{Caffarelli-Silvestre}
have recently given a similar characterization of the Laplacian of
order $\sigma$ in terms of the so-called $\sigma$-harmonic
extension, which is the solution of an elliptic problem with a
degenerate or singular weight. However, even with this local
characterization at hand, many of the proofs that we gave for
$\sigma=1$ cannot be adapted to cope with a general $\sigma$. Hence
we have needed to use new tools, which in several cases do not
involve the extension.  These techniques, which include some new functional inequalities, have also allowed us to improve the results obtained in \cite{dePablo-Quiros-Rodriguez-Vazquez}
 for the case $\sigma=1$.  The case of a bounded domain was not
 treated in that paper.

\section{Main results}\label{sect-main-results}
\setcounter{equation}{0}

As in the case of the PME, there is a unified theory of existence
and uniqueness of a suitable concept of \emph{weak} solution above a critical exponent, given for a general $\sigma\in(0,2)$ by
$m_*\equiv(N-\sigma)_+/N$. The linear case $m=1$ lies always in this range.

\begin{figure}[ht]
\begin{center}
\begin{tabular}{ccc}
\psfig{file=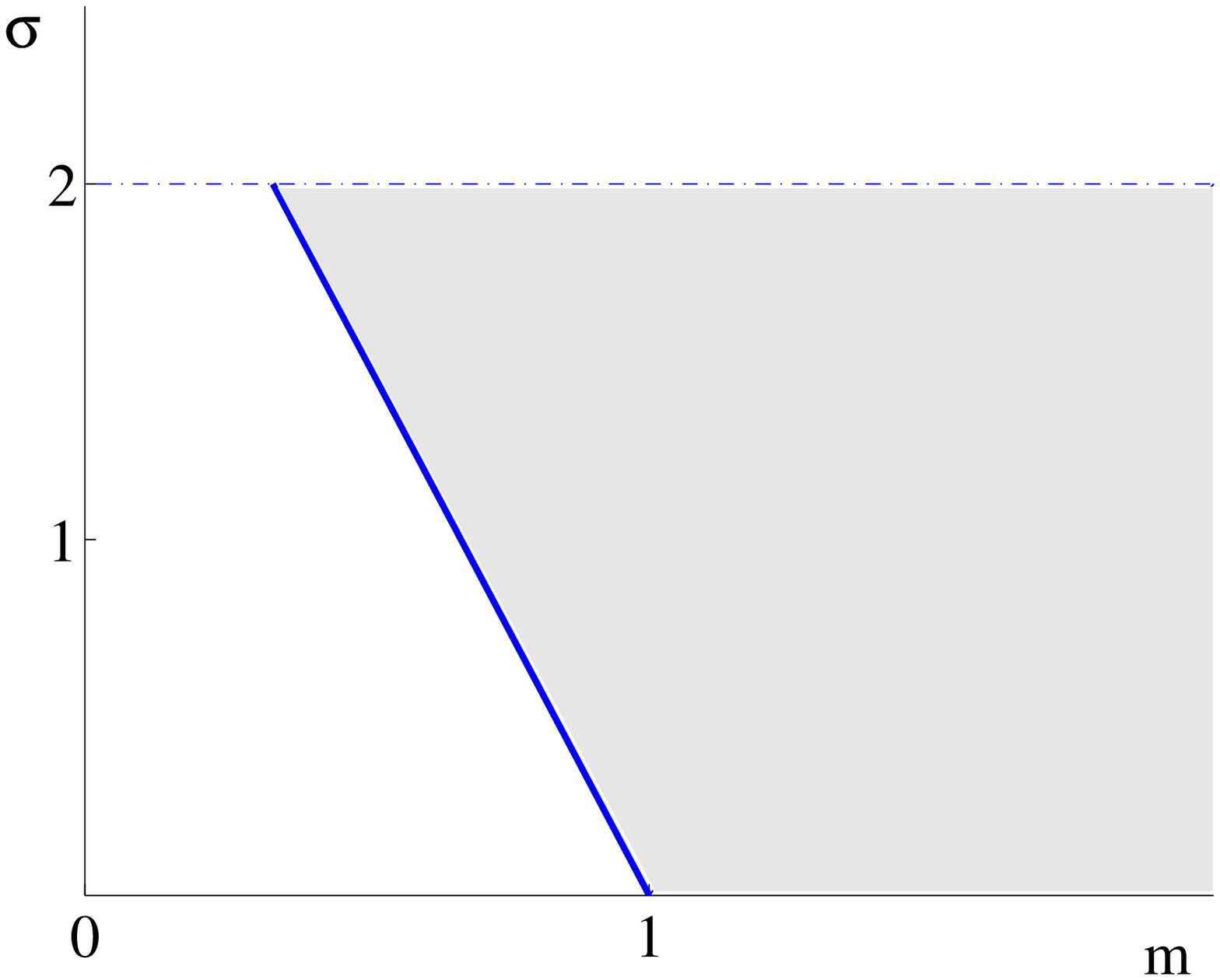,width=4.5cm}&
\psfig{file=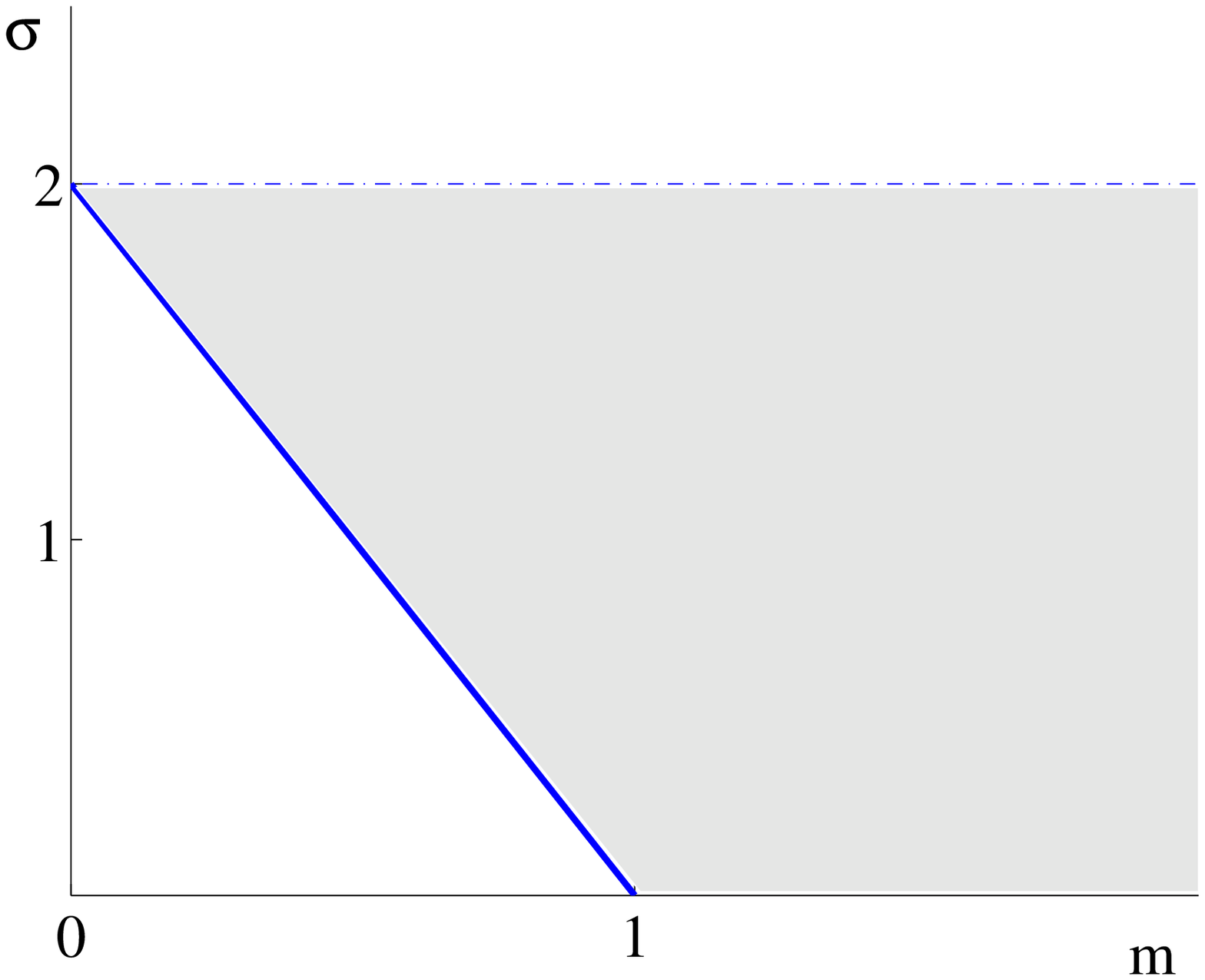,width=4.5cm}&
\psfig{file=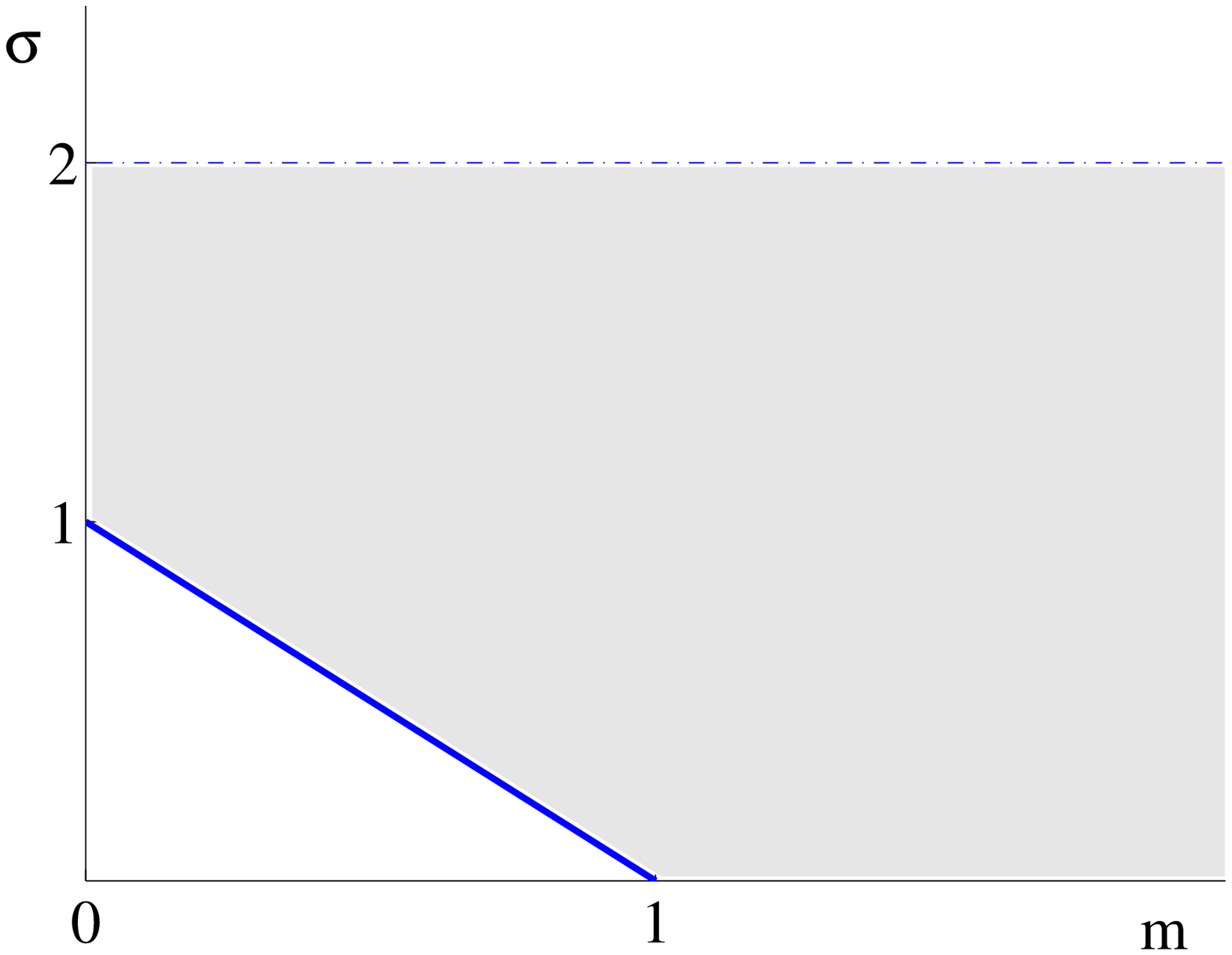,width=4.5cm}
\end{tabular}
\end{center}
\caption{The critical line $m_*=(N-\sigma)/N$ and the supercritical
region $m>m_*$ for $N\ge3$, $N=2$, and $N=1$.}
\end{figure}

\begin{Theorem}\label{th:general}
Let $m>m_*$ and $\sigma\in (0,2)$. For every $f\in L^1(\mathbb{R}^N)$ there
exists a unique weak solution of Problem~\eqref{eq:main}.
\end{Theorem}

The precise definition of weak solution that guarantees uniqueness is stated in
Definition \ref{def:weak.solution.nonlocal} or equivalently in Definition
\ref{def:weak.solution}.

The construction of the solution in the
previous theorem follows from a double limit procedure, approximating first the
initial datum by a sequence of bounded functions, and also approximating
$\mathbb{R}^N$ by a sequence of bounded domains with null boundary data. In
this respect we show existence of a weak solution to the associated
Cauchy-Dirichlet problem, a result that has an independent interest, see
Section~\ref{sec.bounded.domain}.

The weak solutions to Problem~\eqref{eq:main} have some nice qualitative
properties that are summarized as follows.

\begin{Theorem}\label{th:properties}
Assume the hypotheses of Theorem~{\rm \ref{th:general}}, and let $u$
be the weak solution to Problem~\eqref{eq:main}.
\begin{itemize}
\item[\rm (i)]  $\partial_tu\in L^\infty((\tau,\infty):
L^1(\mathbb{R}^N))$ for every $\tau>0$.

\item[\rm (ii)] Mass is conserved:
$\int_{\mathbb{R}^N}u(x,t)\,dx=\int_{\mathbb{R}^N}f(x)\,dx$ for all
$t\ge0$.

\item[\rm (iii)] Let $u_1,u_2$ be the weak solutions to
Problem~\eqref{eq:main} with initial data $f_1,f_2\in
L^1(\mathbb{R}^N)$. Then,
$$
\int_{\mathbb{R}^N}(u_1-u_2)_+(x,t)\,dx\le
\int_{\mathbb{R}^N}(f_1-f_2)_+(x)\,dx,\
\text{\emph{($L^1$-order-contraction property)}}.
$$

\item[\rm (iv)] Any $L^p$-norm  of the solution, $1\le p\le \infty$, is
non-increasing in time.

\item[\rm (v)] The solution is bounded in
$\mathbb{R}^N\times[\tau,\infty)$ for every $\tau>0$
(\emph{$L^1$-$L^\infty$ smoothing effect}). Moreover, for all
$p\ge1$,
\begin{equation}
\|u(\cdot,t)\|_{L^\infty(\mathbb{R}^N)}\le C\,t^{-\gamma_p
}\|f\|_{L^p(\mathbb{R}^N)}^{\delta_p}, \label{eq:L-inf-p}
\end{equation}
with $\gamma_p=(m-1+{\sigma}p /N)^{-1}$, $\delta_p=\sigma
p\gamma_p/N$, and \ $C=C(m,p,N,\sigma)$.

\item[\rm (vi)] If $f\ge0$ the solution is positive  for all $x$ and all
positive times.

\item[\rm (vii)] If either  $m\ge 1$ or $f\ge0$,  then $u\in
C^\alpha(\mathbb{R}^N\times(0,\infty))$  for some $0<\alpha<1$.

\item[\rm (viii)] The solution depends continuously on the
parameters $\sigma\in (0,2)$, $m>m_*$, and   $f\in
L^1(\mathbb{R}^N)$  in the norm of the space
$C([0,\infty):\,L^1(\mathbb{R}^N))$.
\end{itemize}
\end{Theorem}

\noindent\emph{Remarks. } (a)  Properties (i) and (ii) were only known for $\sigma=1$ in the case of nonnegative initial data \cite{dePablo-Quiros-Rodriguez-Vazquez}.

\smallskip

\noindent (b) The positivity property (vi) is not true for the PME in the range $m>1$. The fact that it holds for Problem~\eqref{eq:main} stems from the nonlocal character of the diffusion operator.

\smallskip

\noindent (c) A weak solution satisfying property (i) is said to be a \emph{strong} solution, see
Definition~\ref{def:strong.solution}. These kind of solutions satisfy the equation in
\eqref{eq:main} almost everywhere in $Q=\mathbb{R}^N\times
(0,\infty)$.

Our main interest in this paper is in describing the theory in the above
mentioned range $m>m_*$.  However, we also consider the lower
range $0<m\le m_*$ for contrast. In that range  (which implies that
$0<\sigma<1$ if $N=1$) we obtain existence if we
restrict the data class (or if we relax the concept of solution). In addition,
in order to have uniqueness, we need to ask the solution to be strong.

\begin{Theorem}\label{th:general2} Let $\sigma\in (0,2)$ and $0<m\le m_*$.  For every  $f\in L^1(\mathbb{R}^N)\cap
L^p(\mathbb{R}^N)$  with  $p>p_*(m)=(1-m)N/\sigma$
 there exists a  unique strong solution
 to Problem~\eqref{eq:main}.
\end{Theorem}

\begin{figure}[ht]
\begin{center}
\begin{tabular}{c}
\psfig{file=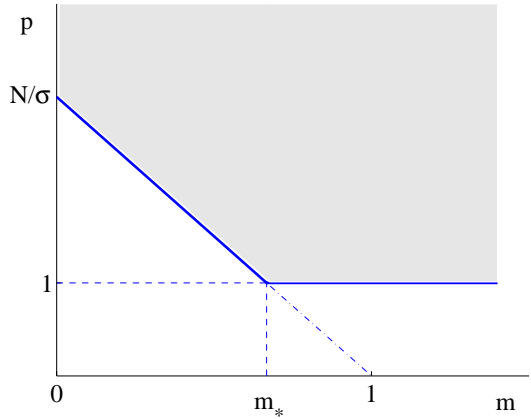,width=7cm}
\end{tabular}
\end{center}
\caption{The critical line $p_*=(1-m)N/\sigma$ and the existence
region of strong solutions.} 
\end{figure}

This theorem improves the results obtained of
\cite{dePablo-Quiros-Rodriguez-Vazquez} for $\sigma=1$, since in
that paper existence and uniqueness of a (weak) solution were only
proved for the case of integrable and bounded initial data, $f\in
L^1(\mathbb{R}^N)\cap L^\infty(\mathbb{R}^N)$. Moreover, the weak
solution was only proved to be strong in the case of nonnegative
initial data.

Some of the properties of the solutions in this lower range are rather
different from those in the upper range.

\begin{Theorem}\label{th:properties2}
Assume the hypotheses of Theorem~{\rm \ref{th:general2}}, and let
$u$ be the strong solution to Problem~\eqref{eq:main}.

\begin{itemize}
\item[\rm (i)] The mass $\int_{\mathbb{R}^N}u(x,t)\,dx$ is conserved if
$m=m_*$.  Mass is not conserved if $m<m_*$. Actually, when $0<m<
m_*$ there is a finite time $T>0$ such that  $u(x,T)\equiv0$ in $\mathbb{R}^N$.

\item[\rm (ii)] There is an $L^1$-order-contraction property.

\item[\rm (iii)] Any $L^p$-norm  of the solution, $1\le p\le \infty$, is
non-increasing in time.

\item[\rm (iv)] The solution is bounded in
$\mathbb{R}^N\times[\tau,\infty)$ for every $\tau>0$
\emph{($L^p$-$L^\infty$ smoothing effect)}. Moreover, if $p>p_*(m)$
(which is necessary to make $\gamma_p>0$) then
formula~\eqref{eq:L-inf-p} holds.

\item[\rm (v)] If $f\ge0$ the solution is positive  for all $x$ and all
positive times up to the extinction time.

\item[\rm (vi)] Let $f\ge0$ and let $T$ be the extinction time. Then
$u\in C^\alpha(\mathbb{R}^N\times(0,T))$ for some $0<\alpha<1$.
\end{itemize}
\end{Theorem}

We will make a series of comments on these two results. First,  we point out
that the result on the conservation of mass for $m=m_*$ is new even for
$\sigma=1$.

A more essential observation is that there is an alternative approach:
 using the results of  Crandall and Pierre \cite{Crandall-Pierre}
we can obtain the existence of a unique so-called \emph{mild} solution for all
$f\in L^1(\mathbb{R}^N)$ in the  whole range of $m$ and $\sigma$, via the
abstract theory of accretive operators. This approach has therefore  the
advantage of having general scope. Two problems arise with this way of looking at the problem:  (a)~how to characterize the mild solution in differential terms; and (b)~how to
derive its properties. Our paper answers both
questions.

The strong solutions that we have constructed are mild solutions. Hence, in our
restricted range of initial data the unique mild solution is a strong solution.
For a general $f\in L^1(\mathbb{R}^N)$ and $m\le m_*$, we will show that the
mild solution is  a \emph{very weak} solution (a solution in the sense of
distributions), see Theorem~\ref{th:mild.very.weak}. However,  we will fail to
prove that this very weak solution is a weak solution  in the sense of
Definition \ref{def:weak.solution.nonlocal}, and hence we will not be able to
obtain the properties listed above (bit note that the $L^1$-contraction holds
since it is a consequence of the accretivity of the operator).

As to the continuous dependence of the solution in terms of the
parameters, convergence in $L^1(\mathbb{R}^N)$ is not expected to
hold for $0<m<m_*$, since mass is not conserved in that region.
Instead we expect to have continuity in weighted spaces, much in the
spirit of~\cite{Benilan-Crandall}. Nevertheless, we are able to
extend the continuous dependence result of
Theorem~\ref{th:properties}-(viii) to cover the case $m=m_*$ for
$N>2$, see Proposition~\ref{th:cont}. We also show that
continuity holds in the upper limit $\sigma\to2$, thus recovering
the standard PME, see Theorem~\ref{th:cont.sigma2}.

\section{Weak solutions. An equivalent problem}\label{sec.prelim}
\setcounter{equation}{0}

\subsection{Weak solutions}

If $\psi$ and $\varphi$ belong to the Schwartz class, Definition
\eqref{def-fourier} of the fractional Laplacian together with
Plancherel's theorem yield
\begin{equation}\label{eq:integration.by.parts}
\int_{\mathbb{R}^N}(-\Delta)^{\sigma/2}\psi\,\varphi=\int_{\mathbb{R}^N}|\xi|^\sigma
  \hat\psi\,\hat\varphi=\int_{\mathbb{R}^N}|\xi|^{\sigma/2} \hat
\psi|\xi|^{\sigma/2}\,\hat
\varphi=\int_{\mathbb{R}^N}(-\Delta)^{\sigma/4}\psi\,(-\Delta)^{\sigma/4}\varphi.
\end{equation}
Therefore, if we multiply  the equation in \eqref{eq:main} by a test function
$\varphi$ and integrate by parts, as usual, we  obtain
\begin{equation}\label{weak-nonlocal}
\displaystyle \int_0^\infty\int_{\mathbb{R}^N}u\dfrac{\partial
\varphi}{\partial
t}\,dxds-\int_0^\infty\int_{\mathbb{R}^N}(-\Delta)^{\sigma/4}(|u|^{m-1}u)(-\Delta)^{\sigma/4}\varphi\,d
xds=0.
\end{equation}
This identity will be the  basis of our definition of a weak solution.

The integrals in \eqref{weak-nonlocal} make sense if $u$ and $|u|^{m-1}u$
belong to suitable spaces. The  correct space  for $|u|^{m-1}u$ is the
fractional Sobolev space $\dot{H}^{\sigma/2}(\mathbb{R}^N)$, defined as the
completion of $C_0^\infty(\mathbb{R}^N)$ with the norm
$$
  \|\psi\|_{\dot{H}^{\sigma/2}}=\left(\int_{\mathbb{R}^N}
|\xi|^\sigma|\hat{\psi}|^2\,d\xi\right)^{1/2}
=\|(-\Delta)^{\sigma/4}\psi\|_{2}.
$$

\begin{Definition}\label{def:weak.solution.nonlocal} A function $u$
is a {\it weak} ($L^1$-energy solution) to Problem \eqref{eq:main}
if:
\begin{itemize}
\item $u\in C([0,\infty): L^1(\mathbb{R}^N))$, $|u|^{m-1}u \in L^2_{\rm
loc}((0,\infty):\dot{H}^{\sigma/2}(\mathbb{R}^N))$;
\item  identity  \eqref{weak-nonlocal}
holds for every $\varphi\in C_0^1(\mathbb{R}^N\times(0,\infty))$;
\item  $u(\cdot,0)=f$ almost everywhere.
\end{itemize}
\end{Definition}

For brevity we will call weak solutions the solutions obtained below according
to this definition, but the complete name {\sl weak $L^1$-energy solution} is
used in the statement to recall that we are making a very definite choice.

The main disadvantage in using this definition is that there is no
formula for the fractional Laplacian of a product, or of a
composition of functions. Moreover,  there is no benefit
 in using compactly supported test functions since their
fractional Laplacian loses this property. To overcome these
difficulties, we will use the fact that our solution $u$ is the
trace of the solution of a \emph{local} problem obtained by
extending $|u|^{m-1}u$ to a half-space whose boundary is our
original space.

\subsection{Extension Method}

 If $g=g(x)$ is a smooth bounded
function defined in $\mathbb{R}^N$, its $\sigma$-harmonic extension
to the upper half-space, $v=\E(g)$, is  the unique smooth bounded
solution $v=v(x,y)$ to
\begin{equation}
\left\{
\begin{array}{ll}
\nabla\cdot(y^{1-\sigma}\nabla v)=0,\qquad &\text{in
}\mathbb{R}^{N+1}_+\equiv\{(x,y)\in\mathbb{R}^{N+1}:
x\in\mathbb{R}^N,
y>0\},\\
v(x,0)=g(x),\qquad&x\in\mathbb{R}^N.
\end{array}
\right. \label{sigma-extension}
\end{equation}
Then, Caffarelli and Silvestre \cite{Caffarelli-Silvestre} proved
that
\begin{equation}
\label{fract-lapla}
-\mu_{\sigma}\lim_{y\to0^+}y^{1-\sigma}\frac{\partial v}{\partial
y}=(-\Delta)^{\sigma/2} g(x),
\quad \mu_{\sigma}={2^{\sigma-1}\Gamma(\sigma/2)}/{\Gamma(1-\sigma/2)}.
\end{equation}
In \eqref{sigma-extension} the
operator  $\nabla$ acts in all $(x,y)$ variables, while in
\eqref{fract-lapla} $(-\Delta)^{\sigma/2}$ acts only on the
$x=(x_1,\cdots,x_N)$ variables.  In the sequel we denote
$$
  L_\sigma v\equiv \nabla\cdot(y^{1-\sigma}\nabla v),\qquad
  \dfrac{\partial v}{\partial y^\sigma}\equiv
  \mu_{\sigma}\lim_{y\to0^+}y^{1-\sigma}\frac{\partial v}{\partial
  y}.
$$

Operators like $L_\sigma$, with a coefficient $y^{1-\sigma}$, which
belongs to the Muckenhoupt space of weights $A_2$ if $0<\sigma<2$,
have been studied by Fabes et al.~in \cite{Fabes-Kenig-Serapioni}.
We make use of this theory later in the proof of positivity, see
Theorem~\ref{th:positivity}.

\medskip

With the above in mind, we rewrite Problem~\eqref{eq:main}  as a
quasi-stationary problem for $w=\E(|u|^{m-1}u)$ with a dynamical
boundary condition
\begin{equation}
\left\{
\begin{array}{ll}
L_\sigma w=0,\qquad &(x,y)\in\mathbb{R}^{N+1}_+,\, t>0,\\
\dfrac{\partial w}{\partial y^\sigma}-\dfrac{\partial
|w|^{\frac1m-1}w}{\partial
t}=0,\qquad&x\in\mathbb{R}^{N},y=0,\, t>0,\\
w=|f|^{m-1}f,\qquad&x\in\mathbb{R}^{N}, y=0,\, t=0.
\end{array}
\right. \label{pp:local}
\end{equation}
To define a weak solution of this problem we multiply formally the
equation in \eqref{pp:local} by a test function $\varphi$ and
integrate by parts to obtain
\begin{equation}\label{weak-local}
\displaystyle \int_0^\infty\int_{\mathbb{R}^{N}}u\dfrac{\partial
\varphi}{\partial
t}\,dxds-\mu_\sigma\int_0^\infty\int_{\mathbb{R}^{N+1}_+}y^{1-\sigma}\langle\nabla
w,\nabla \varphi\rangle\,dxdyds=0,
\end{equation}
where $u=|\Tr(w)|^{\frac1m-1}\Tr(w)$. This holds on the condition that $\varphi$ vanishes for
$t=0$ and also for large $|x|$, $y$ and $t$. We then introduce
the energy space $X^\sigma(\mathbb{R}^{N+1}_+)$, the completion of
$C_0^\infty(\mathbb{R}^{N+1}_+)$ with the norm
\begin{equation}
  \|v\|_{X^\sigma}=\left(\mu_\sigma\int_{\mathbb{R}^{N+1}_+} y^{1-\sigma}|\nabla
v|^2\,dxdy\right)^{1/2}. \label{norma2}
\end{equation}
The trace operator is well defined in this space, see below.

\begin{Definition}\label{def:weak.solution}
A pair of functions $(u,w)$ is a {\it weak solution} to Problem
\eqref{pp:local} if:
\begin{itemize}
\item $u=|\Tr(w)|^{\frac1m-1}\Tr(w)\in
C([0,\infty):L^{1}(\mathbb{R}^{N}))$, $w\in L^2_{\rm
loc}((0,\infty):X^\sigma(\mathbb{R}^{N+1}_+))$;
\item  identity  \eqref{weak-local}
holds for every $\varphi\in C_0^1\left(\overline{\mathbb{R}^{N+1}_+}\times(0,\infty)\right)$;
\item  $u(\cdot,0)=f$ almost everywhere.
\end{itemize}
\end{Definition}

For brevity we will refer sometimes to the solution as only $u$, or
even only $w$, when no confusion arises, since it is clear how to
complete the pair from one of the components, $u=|\Tr(w)|^{\frac1m-1}\Tr(w)$,
$w=\E(|u|^{m-1}u)$.

\medskip

The extension operator is well defined in $\dot{H}^{\sigma/2}(\mathbb{R}^{N})$. It has
an explicit expression using a $\sigma$-Poisson kernel, and
$\E:\dot{H}^{\sigma/2}(\mathbb{R}^{N})\to X^\sigma(\mathbb{R}^{N+1}_+)$ is an isometry, see
\cite{Caffarelli-Silvestre}. The trace operator, $\Tr:
X^\sigma(\mathbb{R}^{N+1}_+)\to\dot{H}^{\sigma/2}(\mathbb{R}^{N})$ is surjective and continuous.
Actually, we have the trace embedding
\begin{equation}
\label{eq:trace.embedding}
 \|\Phi\|_{X^\sigma}\ge\|\E(\Tr(\Phi))\|_{X^\sigma}=\|\Tr(\Phi)\|_{\dot H^{\sigma/2}}
\end{equation}
for any $\Phi\in X^\sigma(\mathbb{R}^{N})$.

\subsection{Equivalence of the weak formulations}

The key point of the above discussion is that the definitions of
weak solution for our original nonlocal  problem and for the
extended local problem are equivalent. Thus, in the sequel we will
switch from one formulation to the other whenever this may offer
some advantage.

\begin{Proposition}\label{th:equivalence}
A function  $u$ is a weak solution to Problem \eqref{eq:main}  if and only if
$(u,\E(|u|^{m-1}u))$ is a weak solution to Problem~\eqref{pp:local}.
\end{Proposition}

Since  $\E:\dot{H}^{\sigma/2}(\mathbb{R}^{N})\to X^\sigma(\mathbb{R}^{N+1}_+)$ is an isometry, we
have
$$
\mu_\sigma\int_{\mathbb{R}^{N+1}_+} y^{1-\sigma}\langle\nabla
\E(\psi),\nabla
  \E(\varphi)\rangle=\int_{\mathbb{R}^{N}}(-\Delta)^{\sigma/4}\psi\,(-\Delta)^{\sigma/4}\varphi,
$$
for every $\psi,\,\varphi\in \dot H^{\sigma/2}(\mathbb{R}^{N})$. Hence the
result follows immediately from the next lemma, which states that any
$\sigma$-harmonic function is  orthogonal in $X^\sigma(\mathbb{R}^{N+1}_+)$ to
every function with trace 0 on $\mathbb{R}^{N}$.

\begin{Lemma}\label{lem:equivalence}
Let $\psi\in\dot H^{\sigma/2}(\mathbb{R}^{N})$ and $\Phi_1,\Phi_2\in
X^\sigma(\mathbb{R}^{N+1}_+)$ such that $\Tr(\Phi_1)=\Tr(\Phi_2)$. Then
$$
\mu_\sigma\int_{\mathbb{R}^{N+1}_+} y^{1-\sigma}\langle\nabla \E(\psi),\nabla
  \Phi_1\rangle=
\mu_\sigma\int_{\mathbb{R}^{N+1}_+} y^{1-\sigma}\langle\nabla \E(\psi),\nabla
  \Phi_2\rangle.
$$
\end{Lemma}
\noindent{\it Proof. } Let $h=\Phi_1-\Phi_2$. Since $\E(\psi)$ is
smooth for $y>0$, given $\varepsilon>0$ we have, after integrating
by parts,
$$
\mu_\sigma\int_\varepsilon^\infty\int_{\mathbb{R}^N}
y^{1-\sigma}\langle\nabla \E(\psi),\nabla
h\rangle\,dxdy=\mu_\sigma\int_{\mathbb{R}^N}
\varepsilon^{1-\sigma}\dfrac{\partial \E(\psi)}{\partial
y}(x,\varepsilon)h(x,\varepsilon)\,dx.
$$
The left-hand side converges to $ \mu_\sigma\int_{\mathbb{R}^{N+1}_+}
y^{1-\sigma}\langle\nabla \E(\psi),\nabla h\rangle$, while the right
hand side tends to 0, since identity \eqref{fract-lapla} holds in
the weak sense in  $H^{-\sigma/2}(\mathbb{R}^N)$, and
$\Tr(h)=0$.~\qed

\section{The problem in a bounded
domain}\label{sec.bounded.domain}
\setcounter{equation}{0}

As an intermediate step in the development of the theory for
Problem~\eqref{eq:main}, we will also consider the Cauchy-Dirichlet problem associated to
the fractional PME,
\begin{equation}
\label{eq:main.bounded} \left\{
\begin{array}{ll}
\dfrac{\partial u}{\partial t} + (-\Delta)^{\sigma/2} (|u|^{m-1}u)=0, &
\qquad x\in\Omega,\; t>0,
\\ [4mm]
u=0, & \qquad  x\in\partial\Omega,\; t>0,\\ [4mm]
u(x,0) = f(x), & \qquad x\in\Omega,%
\end{array}
\right.
\end{equation}
where $\Omega\subset\mathbb{R}^N$ is a smooth bounded domain, $f\in
L^1(\Omega)$.  This problem has an interest in itself.

Let us present here the main facts and results about this problem. In view of
formula \eqref{laplace}, the fractional operator $(-\Delta)^{\sigma/2}$ in a
bounded domain can be described in terms of a spectral decomposition. Let
$\{\varphi_k,\lambda_k\}_{k=1}^\infty$ denote an orthonormal basis of
$L^2(\Omega)$ consisting of eigenfunctions of $-\Delta$ in $\Omega$ with
homogeneous Dirichlet boundary conditions and their corresponding eigenvalues.
The operator $(-\Delta)^{\sigma/2}$ is defined for any $u\in
C_0^\infty(\Omega)$, $ u=\sum_{k=1}^\infty u_k\varphi_k$, by
$$
(-\Delta)^{\sigma/2} u=\sum_{k=1}^\infty\lambda_k^{\sigma/2} u_k\varphi_k.
$$
This operator can be extended by density for $u$ in the Hilbert space
$$
H_0^{\sigma/2}(\Omega)=\{u\in L^2(\Omega):
\|u\|^2_{H_0^{\sigma/2}}=\sum_{k=1}^\infty\lambda_k^{\sigma/2}u_k^2<\infty\}.
$$

\begin{Definition}\label{def:weak.solution.nonlocal.bounded}
A function $u$ is a {\it weak solution} to Problem
\eqref{eq:main.bounded} if:
\begin{itemize}
\item $u\in
C([0,\infty):L^{1}(\Omega))$, $|u|^{m-1}u \in L^2_{\rm
loc}((0,\infty):H_0^{\sigma/2}(\Omega))$;

\item  Identity
$$
\displaystyle
\int_0^T\int_{\Omega}u\dfrac{\partial \varphi}{\partial
t}\,dxds-\int_0^T\int_{\Omega}(-\Delta)^{\sigma/4}u^m(-\Delta)^{\sigma/4}\varphi\,d
xds=0
$$
holds for every $\varphi\in C_0^1(\Omega\times(0,T))$;
\item  $u(\cdot,0)=f$ almost everywhere in $\Omega$.
\end{itemize}
\end{Definition}

The hypotheses that we need in order to have existence when the spatial domain
is bounded coincide with the ones we have when the spatial domain is the whole
$\mathbb{R}^N$.

\begin{Theorem}\label{th:existence-bounded}
Problem~\eqref{eq:main.bounded} has a unique weak solution if
$m>m_*$ and $f\in L^1(\Omega)$, which is moreover strong, and a
unique strong solution if $m\le m_*$ and $f\in L^p(\Omega)$ with
$p>p_*(m)=(1-m)N/\sigma$.
\end{Theorem}
As for the properties of the solutions, most of them, though not all, coincide
with the ones that hold when the domain is the whole space.
\begin{Theorem}
\label{th:properties-bounded} Assume the hypotheses of
Theorem~\ref{th:existence-bounded}, and let $u$ be the strong solution to
Problem~\eqref{eq:main.bounded}.

\begin{enumerate}
\item[\rm (i)] The solution is bounded in
$\Omega\times[\tau,\infty)$ for every $\tau>0$. Moreover, a formula
analogous to~\eqref{eq:L-inf-p} holds.

\item[\rm (ii)] As a consequence,  $ \int_{\Omega}u(x,t)\,dx=O(t^{-\gamma_p})$.
Moreover, if $0<m<1$ there is extinction in finite time.

\item[\rm (iii)] There is an $L^1$-order-contraction property.

\item[\rm (iv)] Any $L^p$-norm  of the solution, $1\le p\le \infty$, is
non-increasing in time.

\item[\rm (v)] The solution depends continuously on the parameters
$\sigma\in (0,2)$, $m>m_*$, and   $f\in L^1(\Omega)$  in the norm of
the space $C([0,\infty):\,L^1(\Omega))$.
\end{enumerate}
\end{Theorem}

The results of \cite{Athanasopoulos-Caffarelli} imply that $u\in
C^\alpha$ for $m\ge1$. Positivity for any $m>0$ when the initial
data are nonnegative, and $C^\alpha$ regularity for $m<1$ are still
open problems.

The construction of a solution uses the analogous to the
Caffarelli-Silvestre extension~\eqref{sigma-extension}, restricted
here to the upper half-cylinder $C_\Omega=\Omega\times(0,\infty)$,
with null condition on the lateral boundary, $\partial
\Omega\times(0,\infty)$,  a construction considered in
\cite{Cabre-Tan}, \cite{Stinga-Torrea},
\cite{Brandle-Colorado-dePablo-Sanchez}, \cite{capella-d-d-s}. Thus,
$w=\E(|u|^{m-1}u)$ satisfies
\begin{equation}
\label{pp:local.bounded} \left\{
\begin{array}{ll}
L_\sigma w=0,\qquad &(x,y)\in C_\Omega,\, t>0,\\
w=0,\qquad&x\in\partial\Omega, y>0,\,t>0,\\
 \dfrac{\partial
w}{\partial y^\sigma}-\dfrac{\partial |w|^{\frac1m-1}w}{\partial
t}=0,\qquad&x\in\Omega, y=0,\, t>0,\\
w=|f|^{m-1}f,\qquad&x\in\Omega, y=0,\, t=0.
\end{array}
\right.
\end{equation}

In order to define a weak solution to \eqref{pp:local.bounded} we
need to consider the space $X^\sigma_0({C_\Omega})$, the closure of
$C^\infty_c(C_\Omega)$ with respect to the norm \eqref{norma2}, with
$\mathbb{R}^{N+1}_+$ substituted by $C_\Omega$. The extension and
trace operators between $H_0^{\sigma/2}(\Omega)$ and
$X^\sigma_0({C_\Omega})$ satisfy the same properties as in the case
of the whole space. In fact
$$
\begin{array}{rl}
\displaystyle\left(\mu_\sigma\int_{C_\Omega} y^{1-\sigma}|\nabla
\E(\varphi)|^2\right)^{1/2}&= \displaystyle
\|\E(\varphi)\|_{X^\sigma_0}=\displaystyle\|\varphi\|_{H_0^{\sigma/2}}
\\ [3mm] &= \|(-\Delta)^{\sigma/4}
\varphi\|_{2}=\left(\sum_{k=1}^\infty\lambda_k^{\sigma/2}
\varphi_k^2\right)^{1/2}\,.
\end{array}$$
See for instance \cite{Brandle-Colorado-dePablo-Sanchez} for the
explicit expression of $\E(\varphi)$ in terms of the coefficients
$\varphi_k$.

\begin{Definition}\label{def:weak.solution.local.bounded}
A pair of functions $(u,w)$ is a {\it weak solution} to Problem
\eqref{pp:local.bounded} if:
\begin{itemize}
\item $u=|\Tr(w)|^{\frac1m-1}\Tr(w)\in
C([0,\infty):L^{1}(\Omega))$, $w\in L^2_{\rm
loc}((0,\infty):X_0^\sigma(C_\Omega))$;
\item  identity
$$
\displaystyle \int_0^\infty\int_{\Omega}u\dfrac{\partial
\varphi}{\partial
t}\,dxds-\mu_\sigma\int_0^\infty\int_{C_\Omega}y^{1-\sigma}\langle\nabla
w,\nabla \varphi\rangle\,dxdyds=0,
$$
holds for every $\varphi=\varphi(x,y,t)$ such that $\varphi\in
C_0^1\left(\Omega\times[0,\infty)\times(0,\infty)\right)$;
\item  $u(\cdot,0)=f$ almost everywhere in $\Omega$.
\end{itemize}
\end{Definition}

As it happens for the case where $\Omega=\mathbb{R}^N$, the two
definitions of weak solution,
Definitions~\ref{def:weak.solution.nonlocal.bounded} and
\ref{def:weak.solution.local.bounded},  are equivalent.

\medskip

\noindent{\it Remark.} The space $H_0^{\sigma/2}(\Omega)$ can also
be defined by interpolation, see \cite{Lions-Magenes}.  We notice
that, though for $1<\sigma<2$ the solutions are zero almost
everywhere at the boundary, for $0<\sigma\le1$ the functions in
$H_0^{\sigma/2}(\Omega)$ do not have a trace, \cite{Lions-Magenes}.
Therefore, the boundary condition must be understood in a weak
sense, see also \cite{Cabre-Tan} for the case $\sigma=1$.

\section{Some functional inequalities}\label{sec.functional.inequalities}
\setcounter{equation}{0}

In this section we gather some functional inequalities related with the
fractional Laplacian, both in the whole space or in a bounded
domain,  that will play an important role in the sequel. The first
one, Strook-Varopoulos' inequality, is well known in the theory of
sub-Markovian operators \cite{LS}. Nevertheless, we give a very short proof
using the extension operator that makes apparent the power of this technique.

\begin{Lemma}[Strook-Varopoulos' inequality] Let $0<\gamma<2$, $q>1$.
Then
\begin{equation}\label{eq:strook.varopoulos}
\int_{\mathbb{R}^N}(|v|^{q-2}v)(-\Delta)^{\gamma/2} v \ge
\frac{4(q-1)}{q^2}\int_{\mathbb{R}^N}\left|(-\Delta)^{\gamma/4}|v|^{q/2}\right|^2
\end{equation}
for all $v\in L^q(\mathbb{R}^N)$ such that $(-\Delta)^{\gamma/2}v\in
L^q(\mathbb{R}^N) $.
\label{lem:S-V}
\end{Lemma}

\noindent{\it Proof. } Using property
\eqref{eq:integration.by.parts} and Lemma~\ref{lem:equivalence}, we
get
$$
\begin{array}{rcl}
\displaystyle\int_{\mathbb{R}^N}(|v|^{q-2}v)(-\Delta)^{\gamma/2} v
&=&
\displaystyle\int_{\mathbb{R}^N}(-\Delta)^{\gamma/4}(|v|^{q-2}v)(-\Delta)^{\gamma/4}
v\\ [3mm]
&=&\displaystyle\mu_\gamma\int_{\mathbb{R}^{N+1}_+}y^{1-\gamma}\langle
\nabla(|\E(v)|^{q-2}\E(v)),\nabla \E(v)\rangle\\ [3mm]
&=&\displaystyle\mu_\gamma
\frac{4(q-1)}{q^2}\int_{\mathbb{R}^{N+1}_+}y^{1-\gamma}|\nabla
(|\E(v)|^{q/2})|^2 \\ [3mm] &\ge&\displaystyle
\frac{4(q-1)}{q^2}\int_{\mathbb{R}^N}\left|(-\Delta)^{\gamma/4}|v|^{q/2}\right|^2.
\end{array}
$$
In the last step we get only inequality because the function
$|\E(v)|^{q/2}$ is not necessarily $\gamma$-harmonic. \qed

With the same technique a generalization of \eqref{eq:strook.varopoulos} can be
proved.

\begin{Lemma} Let $0<\gamma<2$. Then
\begin{equation}\label{eq:strook.varopoulos2}
\int_{\mathbb{R}^N}\psi(v)(-\Delta)^{\gamma/2} v \ge
\int_{\mathbb{R}^N}\left|(-\Delta)^{\gamma/4}\Psi(v)\right|^2
\end{equation}
whenever $\psi'=(\Psi')^2$.
\label{lem:S-V2}\end{Lemma}

\noindent{\it Proof. }
Use the extension method and the property $\langle\nabla\psi(w),\nabla
w\rangle=|\nabla\Psi(w)|^2$.
\qed

In order to prove our second inequality we need the well-known
Hardy-Littlewood-Sobolev's inequality \cite{Hardy-Littlewood},
\cite{Sobolev}, \cite{Lieb}: for every $v$ such that
$(-\Delta)^{\gamma/2}v\in L^r(\mathbb{R}^N)$, $1<r<N/\gamma$, $0<\gamma<2$, it
holds
\begin{equation}
\|v\|_{r_1}\le
c(N,r,\gamma)\|(-\Delta)^{\gamma/2}v\|_r, \label{HLS}
\end{equation}
$r_1=\frac{Nr}{N-\gamma r}$. Putting for instance $r=2$, $\gamma=\sigma/2$, we
obtain the inclusion $\dot{H}^{\sigma/2}(\mathbb{R}^N)\hookrightarrow
L^{\frac{2N}{N-\sigma}}(\mathbb{R}^N)$ whenever $N>\sigma$. What happens for
$N=1\le\sigma<2$? Or more generally, for $r\ge N/\gamma$?  We
answer this question in the next lemma.

\begin{Lemma}[Nash-Gagliardo-Nirenberg type inequality]
Let $p\ge1$, $r>1$, $0<\gamma<\min\{N,2\}$. There is a constant
$C=C(p,r,\gamma,N)>0$ such that for any $v\in L^p(\mathbb{R}^N)$ with
$(-\Delta)^{\gamma/2}v\in L^r(\mathbb{R}^N)$ we have
  \begin{equation}
   \label{eq:gagliardo.nirenberg.type.inequality}
    \|v\|_{r_2}^{\alpha+1}\le
    C\|(-\Delta)^{\gamma/2}v\|_r\,\|v\|_p^\alpha\,,
 \end{equation}
 where $r_2=\frac{N(rp+r-p)}{r(N-\gamma)}$, $\alpha=\frac{p(r-1)}r$.
\label{lem-NGN}\end{Lemma} \noindent{\it Proof. } We  use
\eqref{eq:strook.varopoulos}. Estimate now the left hand side of this inequality using
inequality  \eqref{HLS}, and the right hand side with H\"older's
inequality, to get
\eqref{eq:gagliardo.nirenberg.type.inequality}.~\qed

Notice that, for $r=2$ and $\gamma=\sigma/2$, this corresponds to
using inequality \eqref{HLS} for the space
$\dot{H}^{\sigma/4}(\mathbb{R}^N)$ instead of
$\dot{H}^{\sigma/2}(\mathbb{R}^N)$, thus allowing all values of
$\sigma\in(0,2)$ even in the case $N=1$.

Inequalities of this kind are already available \cite{BKM}. However, this
particular one is, up to our knowledge, new.  Let us explain in more detail the consequences of this inequality in relation to \eqref{HLS}.

Assume first that $N>\gamma r$. Hardy-Littlewood-Sobolev's
inequality \eqref{HLS} implies that if $(-\Delta)^{\gamma/2}v\in
L^r$ then $v\in L^{r_1}$. Assuming also $v\in L^p$, then
\eqref{eq:gagliardo.nirenberg.type.inequality} gives that $v\in
L^{r_1}\cap L^{r_2}$, which is always stronger that $v\in
L^{r_1}\cap L^{p}$. Both results coincide in the case $p=r_1=r_2$.

If on the contrary $N\le\gamma r$, we cannot apply Hardy-Littlewood-Sobolev,
but
\eqref{eq:gagliardo.nirenberg.type.inequality} gives that $v\in L^{p}\cap
L^{r_2}$.

\medskip

 We now consider the case of a bounded domain
$\Omega\subset\mathbb{R}^N$.  The characterization of
the fractional Laplacian in terms of the extension to the
half-cylinder $C_\Omega$ allows us to repeat the proofs of
Lemmas~\ref{lem:S-V} and \ref{lem:S-V2} in the case where the domain
is bounded.

\begin{Lemma} Strook-Varopoulos' inequalities \eqref{eq:strook.varopoulos} and \eqref{eq:strook.varopoulos2} hold true with $\mathbb{R}^N$
substituted by $\Omega\subset\mathbb{R}^N$ bounded.
\end{Lemma}

On the other hand, let $v\in H_0^{\sigma/2}(\Omega)$. Consider its
$\sigma$-extension $w=\E(v)$, ($\sigma=2\gamma$), and let
$\widetilde w$ be the extension of $w$ by zero outside the
half-cylinder. Then we have the estimate, see
\cite{Brandle-Colorado-dePablo-Sanchez},
\cite{Cotsiolis-Tavoularis},
$$
\mu_\sigma\int_{\mathbb{R}^{N+1}_+}y^{1-\sigma}|\nabla\widetilde w|^2\ge S(\sigma,N)
\left(\int_{\mathbb{R}^N}|\Tr(\widetilde
w)|^{\frac{2N}{N-\sigma}}\right)^{\frac{N-\sigma}N}.
$$
i.e.,
\begin{equation}
\mu_\sigma\int_{\Omega}y^{1-\sigma}|\nabla w|^2\ge S(\sigma,N)
\left(\int_{\Gamma}|v|^{\frac{2N}{N-\sigma}}\right)^{\frac{N-\sigma}N}.
\label{BCdPS}\end{equation} The left hand side equals
$\|(-\Delta)^{\gamma/2}v\|_2^2$.  That is, we have obtained
inequality \eqref{HLS} in the case $r=2$.

From this point, we can repeat the proof of Lemma~\ref{lem-NGN}, which only
uses the case just proved and H\"older's inequality, thus obtaining inequality
\eqref{eq:gagliardo.nirenberg.type.inequality} also for a bounded domain.

More important is the following application.

\begin{Lemma}[Sobolev type inequality]
\label{lem-NGN2} Let $\Omega\subset\mathbb{R}^N$ be a bounded
domain, and let $v$ be such that $(-\Delta)^{\gamma/2}v\in
L^r(\Omega)$, $N\ge1$, $0<\gamma<2$. Then  we have
  \begin{equation}
   \label{eq:ngn2}
    \|v\|_{q}\le C(q,r,N,\gamma)\|(-\Delta)^{\gamma/2}v\|_r
 \end{equation}
 for every $1\le q\le\frac{Nr}{N-\gamma r}$ if $N>\gamma r$, or for every $q\ge1$ if
 $\gamma<N\le\gamma r$.
\end{Lemma}

\noindent{\it Proof. } If $N>\gamma r$ we just apply
\eqref{eq:gagliardo.nirenberg.type.inequality} with
$p=\frac{Nr}{N-\gamma r}$, and apply H\"older's inequality for the
exponents $1\le q<\frac{Nr}{N-\gamma r}$. If $\gamma<N\le\gamma r$,
and given any $q>\frac N{N-\gamma}$, we
apply~\eqref{eq:gagliardo.nirenberg.type.inequality}, this time with
$s=\frac{Nq}{N+\gamma q}<\frac N\gamma\le r$ instead of $r$, and
$p=q$. We end again with H\"older's inequality for the exponents
$1\le q\le\frac N{N-\gamma}$.~\qed


\section{Uniqueness}\label{sect-uniqueness}
\setcounter{equation}{0}

In this section we prove the uniqueness parts of
Theorems~\ref{th:general} and~\ref{th:general2}.

\noindent\emph{Notations.}  We will use the simplified notation
$u^m$ for data of any sign, instead of the actual  ``odd power''
$|u|^{m-1}u$. In the same way, $w^{1/m}$ will stand for
$|w|^{\frac1m-1}w$. In addition,  $\|\cdot\|_p$ will denote the norm
$\|\cdot\|_{L^p(\mathbb{R}^N)}$  or $\|\cdot\|_{L^p(\Omega)}$.

\subsection{$m>m_*$. Uniqueness of weak solutions}

\begin{Theorem}\label{th:uniqueness1}
Let $f\in L^1(\mathbb{R}^N)$ and $m> m_*$. Problem \eqref{eq:main}
has at most one weak solution.
\end{Theorem}

\noindent{\it Proof. } We adapt the classical uniqueness proof due
to Oleinik \cite{OKC}. This will require $u\in
L^{m+1}({\mathbb{R}^N}\times(0,T))$, which will be true if $m> m_*$.
To prove this we apply H\"older's inequality twice, first in space
and then in time, using inequality \eqref{HLS}. Assume first
$N>\sigma$. We have
$$
\begin{array}{l}
\displaystyle\int_0^T\int_{\mathbb{R}^N} |u|^{m+1}\,dxdt \le
\int_0^T\Big(\int_{\mathbb{R}^N} |u|\,dx\Big)^\beta\,
(\int_{\mathbb{R}^N}
|u|^{\frac{2Nm}{N-\sigma}}\,dx\Big)^{1-\beta}\,dt \\
[3mm]\displaystyle\qquad\le CT^{\mathbb{R}^N}\max_{t\in
[0,T]}\|u(\cdot,t)\|_1^\beta \,
\Big[\int_0^T\Big(\int_{\mathbb{R}^N}
|u|^{\frac{2Nm}{N-\sigma}}\,dx\Big)^{\frac{N-\sigma}{N}}\,dt\Big]^{1-\gamma},
\end{array}
$$
where $\beta=\frac{N(m-1)+\sigma(m+1)}{N(2m-1)+\sigma}$ and
$\gamma=\frac{N(m-1)+\sigma}{N(2m-1)+\sigma}$. Observe that $m>m_*$
implies $\beta,\,\gamma\in(0,1)$. Applying now
inequality~\eqref{HLS}, we get
$$
\int_0^T\int_{\mathbb{R}^N} |u|^{m+1}\,dxdt\le
C\Big[\int_0^T\|u^m(\cdot,t)\|^2_{\dot{H}^{\sigma/2}}\,dt\Big]^{1-\gamma}\le
C.
$$
In the case $N=1$ and $1\le \sigma<2$ the computation is similar.
But we have to use the Nash-Gagliardo-Nirenberg type inequality
\eqref{eq:gagliardo.nirenberg.type.inequality} instead to get the
same conclusion. What we get in this case is
$$
\begin{array}{l}
\displaystyle\int_0^T\int_{\mathbb{R}^N} |u|^{m+1}\,dxdt
\le\displaystyle \left(\int_0^T\Big(\int_{\mathbb{R}^N}
|u|\,dx\Big)^{\frac\delta\gamma}\,dt\right)^{\gamma}
\left(\int_0^T\Big(\int_{\mathbb{R}^N}
|u|^{\frac{2m+1}{2-\sigma}}\,dx\Big)^{2-\sigma}\,dt\right)^{1-\gamma}\\
[3mm]
\qquad\qquad\qquad\quad\qquad\le\displaystyle CT^\gamma\max_{t\in
[0,T]}\|u(\cdot,t)\|_1^{\delta+1-\gamma}\,
\left(\int_0^T\|u^m(\cdot,t)\|^2_{\dot{H}^{\sigma/2}}\,dt\right)^{1-\gamma}\le
C,
\end{array}
$$
where now $\delta=\frac{\sigma(m+1)-1}{2m-1+\sigma}$ and
$\gamma=\frac{m-1+\sigma}{2m-1+\sigma}$.

We now proceed with the core of the proof. Let $u$ and
$\widetilde u$ be two weak solutions to Problem
\eqref{eq:main}. We take the following function as test in the weak
formulation
$$
\varphi(x,t) =\int_t^T (u^m-\widetilde u^m)(
x,s)\,ds,\qquad 0\le t\le T,
$$
with $\varphi\equiv0$ for $t\ge T$.  We have
$$
\begin{array}{l}
\displaystyle\int_0^T\int_{{\mathbb{R}^N}}(u-\widetilde
u)(x,t)(u^m-\widetilde  u^m)(x,t)\,dxdt \\ [4mm] \qquad\quad+
\displaystyle\int_0^T\int_{\mathbb{R}^N}(-\Delta)^{\sigma/4}(u^m-\widetilde
u^m)(x,t)\int_t^T (-\Delta)^{\sigma/4}(u^m-\widetilde u^m)(x,s)
\,ds\,dxdt = 0.
\end{array}
$$
Integration of the second term gives
$$
\begin{array}{l}
\displaystyle\int_0^T\int_{{\mathbb{R}^N}}(u-\widetilde
u)(x,t)(u^m-\widetilde u^m)(x,t)\,dxdt
\\
[4mm] \qquad\qquad\qquad\qquad\qquad\qquad\qquad+\displaystyle\frac
12\int_{{\mathbb{R}^N}}\left(\int_0^T
(-\Delta)^{\sigma/4}(u^m-\widetilde u^m)(x,s)\, ds\right)^2dx= 0.
\end{array}$$
Since both integrands are nonnegative, they must be identically zero. Therefore,
$u=\widetilde u$.~\qed

\noindent{\it Remark.} The same proof works without any restriction
on the exponent $m$ provided $u\in
L^{m+1}({\mathbb{R}^N}\times(\tau,T))$.

\subsection{$m\le m_*$. Uniqueness of strong solutions}

Weak solutions satisfy the equation in~\eqref{eq:main} in the sense
of distributions. Hence, if the left hand side is a function, the
right hand side is also a function and the equation holds almost
everywhere. This fact allows to prove  several important properties,
among them uniqueness for $m\le m_*$, and hence motivates the
following definition.

\begin{Definition} We say that a weak solution $u$  to
Problem \eqref{eq:main} is a strong solution if  $\partial_tu\in
L^\infty((\tau,\infty):L^1(\mathbb{R}^N))$, for every $\tau>0$.
\label{def:strong.solution}\end{Definition}

In the case of strong solutions the uniqueness result also provides
a comparison principle. The following uniqueness proof is valid for
all values of $m>0$.

\begin{Theorem}\label{th:uniqueness2} Let $m>0$. If $u_1$, $u_2$, are strong solutions to Problem~\eqref{eq:main} with initial data
$f_1,f_2\in L^1(\mathbb{R}^N)$, then, for every $0\le t_1<t_2$ it
holds
\begin{equation}
\label{eq:contractivity.uniqueness2}
\int_{{\mathbb{R}^N}}(u_1-u_2)_+(x,t_2)\,dx\le
\int_{{\mathbb{R}^N}}(u_1-u_2)_+(x,t_1)\,dx.
\end{equation}

\end{Theorem}

\noindent{\it Proof. } Let $p\in C^1(\mathbb{R})\cap
L^\infty(\mathbb{R})$ be such that $p(s)=0$ for $s\le0$, $p'(s)>0$ for $s>0$
and $0\le p\le1$, and let $j$ be such that $j'=\sqrt{p'}$, $j(0)=0$. We will
choose $p$ as an approximation to the sign function.

Let us first assume that $t_1>0$. We subtract the equations satisfied by $u_1$ and $u_2$, multiply by the function $\varphi=p(u_1^m-u_2^m)$, and integrate by parts to get
$$
\int_{t_1}^{t_2}\int_{\mathbb{R}^N}\frac{\partial(u_1-u_2)}{\partial
t}p(u_1^m-u_2^m)=-\int_{t_1}^{t_2}\int_{\mathbb{R}^N}
(-\Delta)^{\sigma/2}(u_1^m-u_2^m)p(u_1^m-u_2^m).
$$
We now apply the generalized Strook-Varopoulos inequality
\eqref{eq:strook.varopoulos2}, to get
$$
\int_{t_1}^{t_2}\int_{\mathbb{R}^N}\frac{\partial(u_1-u_2)}{\partial
t}p(u_1^m-u_2^m)\le-C\int_{t_1}^{t_2}\int_{\mathbb{R}^N}
|(-\Delta)^{\sigma/4}j(u_1^m-u_2^m)|^2\le0.
$$
 We end by letting $p$ tend to the sign function. The
case $t_1=0$ is obtained passing to the limit.~\qed

For the problem posed in a bounded domain, the above proofs of uniqueness work
without any change.

\section{Existence for bounded initial data}\label{sect-existence}
\setcounter{equation}{0}

Crandall and Pierre developed in \cite{Crandall-Pierre} an abstract
approach to study evolution equations of the form $\partial_t
u+A\varphi(u)=0$ when $A$ is an $m$-accretive operator in $L^1$ and
$\varphi$ is a monotone increasing real function.  This allows to
obtain a so-called \emph{mild} solution using Crandall-Liggett's
Theorem. Our problem falls within this framework. However, such an
abstract construction does not give enough information to prove that
the mild solution is in fact a weak solution, in other words, to
identify the solutions in a differential sense. We will use an
alternative approach to construct the mild solution whose main
advantage is precisely that it provides enough estimates to show
that it is a weak solution, and later that it is strong.

In order to develop the theory for Problem~\eqref{eq:main}, we will
approximate the initial data $f$ by a sequence $f_n\in
L^1({\mathbb{R}^N})\cap L^\infty({\mathbb{R}^N})$, and use a
contraction property in order to pass to the limit. Hence, our first
task is to obtain existence for integrable, bounded initial data.
This is the goal of the present  section

We will construct solutions by means of Crandall-Liggett's
Theorem~\cite{Crandall-Liggett}, which is based on an implicit in
time discretization. Hence, we will have to deal with the elliptic
problem
\begin{equation}
\left\{
\begin{array}{ll}
L_\sigma w=0\qquad &\mbox{ in } \mathbb{R}^{N+1}_+,\\
-\dfrac{\partial w}{\partial y^\sigma}+w^{1/m}=g\qquad&\mbox{ on }
{\mathbb{R}^N}.
\end{array}
\right. \label{pe:local}
\end{equation}
Equalities on $\mathbb{R}^N$ have to be understood in the sense of
traces. To show existence of a weak solution for this problem we
approximate the domain $\mathbb{R}^{N+1}_+$ by half-cylinders,
$B_R\times(0,\infty)$, with zero data at the lateral boundary,
$\partial B_R\times(0,\infty)$. We recall that in the case
$\sigma=1$ a similar construction is performed in
\cite{dePablo-Quiros-Rodriguez-Vazquez}, though there half-balls
were used instead of half-cylinders.

Now we have two choices: either we first pass to the limit in the
discretization, to obtain a solution of the parabolic problem in the
ball $B_R$, and then pass to the limit $R\to\infty$; or we first
pass to the limit in $R$ to obtain a solution of the elliptic
problem in the whole space and then pass to limit in the
discretization. We will follow both approaches (each of them has its
own advantages) and will later prove that both of them produce the
same solution.

Instead of just considering problems in balls we will analyze the
case of any bounded open domain $\Omega\subset\mathbb{R}^N$, since
it has independent interest.

\subsection{Problem in a bounded domain}\label{ssec.bounded}

In order to check that the hypotheses of Crandall-Liggett's theorem
hold, we have to prove existence of a weak solution $w$  of
\eqref{pe:local} (defined in the standard way) and contractivity of
the map $g\mapsto w^{1/m}(\cdot,0)$ in the norm of $L^1(\Omega)$ for
the elliptic problem for all $g\in L^\infty(\Omega)$.

\begin{Theorem}\label{th:pe}
Let $\Omega\subset\mathbb{R}^N$ be a bounded domain. For every $g\in
L^\infty(\Omega)$ there exists a unique weak solution $w\in
X_0^\sigma(C_\Omega)$ to Problem \eqref{pe:local}. It satisfies
$\|w(\cdot,0)\|_{\infty}\le \|g\|_{\infty}^m$. Moreover, there is a
$T$-contraction property in $L^1$: if $w$ and $\widetilde w$ are the
solutions corresponding to data $g$ and $\widetilde g$, then
\begin{equation}
\int_{\Omega}\left[w^{1/m}(x,0)-\widetilde
w^{1/m}(x,0)\right]_+\,dx\le \int_{\Omega}\left[g(x)-\widetilde
g(x)\right]_+\,dx. \label{eq:contract-e}\end{equation}
\end{Theorem}

\noindent{\it Proof. }  The existence of a weak solution, i.e., a
function $w\in X_0^\sigma(C_\Omega)$ satisfying
\begin{equation}
\label{weak.local.bdd} \mu_\sigma\int_{C_\Omega}
y^{1-\sigma}\langle\nabla
w,\nabla\varphi\rangle+\int_{\Omega}v^{1/m}\varphi-\int_{\Omega}g\varphi=0,
\end{equation}
$v=\Tr(w)$, for every test function $\varphi$, is obtained in a standard way by
minimizing the functional
$$
J(w)=\frac{\mu_\sigma}2\int_{C_\Omega}y^{1-\sigma}|\nabla
w|^2+\frac{m}{m+1}\int_{\Omega}|v|^{\frac{m+1}m}- \int_{\Omega}vg.
$$
This functional is coercive in $X_0^\sigma(C_\Omega)$. Indeed, the
first term $\| w\|^2_{X_0^\sigma}$, and is using H\"older's
inequality, we have
$$
\left|\int_\Omega vg\right|\le
\|v\|_{\frac{2N}{N-\sigma}}\,\|g\|_{\frac{2N}{N+\sigma}}\le
\varepsilon\|v\|^2_{\frac{2N}{N-\sigma}}+\frac1\varepsilon\|g\|_{\frac{2N}{N+\sigma}}^2.
$$
Now, the trace embedding \eqref{BCdPS} implies
\begin{equation}
J(w)\ge C_1 \| w\|^2_{X_0^\sigma}-C_2 .
\label{coercivo}\end{equation}
For $N=1\le\sigma<2$ we use inequality
\eqref{eq:ngn2} instead. In fact,
putting $q=r=2$, $\gamma=\sigma/2$, we  get $\| w\|_{X_0^\sigma} \ge C\|v\|_2$.
We obtain again \eqref{coercivo}.

We now establish contractivity of solutions to Problem
\eqref{pe:local} in $L^1(\Omega)$.  Let $w$ and $\widetilde w$ be
the solutions corresponding to data $g$ and $\widetilde g$. If we
consider in the weak formulation the test function
$\varphi=p(w-\widetilde w)$, where $p$ is any smooth monotone
approximation of the sign function, $0\le p(s)\le 1$, $p'(s)\ge 0$,
we get
$$
\mu_\sigma\int_{C_\Omega}y^{1-\sigma}p'(w-\widetilde
w)|\nabla(w-\widetilde w)|^2+\int_{\Omega}(w^{1/m}-\widetilde
w^{1/m})\,p(w-\widetilde w)-\int_{\Omega}(g-\widetilde
g)\,p(w-\widetilde w)=0.
$$
Passing to the limit, we obtain
$$
\int_{\Omega}(w^{1/m}-\widetilde
w^{1/m})_+\le\int_{\Omega}(g-\widetilde g) \,\mathop{\rm
sign}(w-\widetilde w)\le\int_{\Omega}(g-\widetilde g)_+\,.
$$
In particular, under the assumption $g\ge0$ we have
$w(\cdot,0)\ge0$. Standard comparison gives now $w\ge0$ in
$C_\Omega$. In the same way we can establish a contractivity
property for subsolutions and supersolutions to the problem with
nontrivial, on $\Omega\times(0,\infty)$, boundary condition. We thus
may take the constant function $\widetilde g=\|g\|_\infty$ as a
supersolution, to get $\|w^{1/m}\|_\infty\le \|g\|_\infty$. We also
deduce the estimate
\begin{equation}
  \mu_\sigma\int_{C_\Omega} y^{1-\sigma}|\nabla w|^2\le\int_\Omega gw\le
  C.
\label{energy.bdd}\end{equation} \qed

We now construct a solution to the parabolic problem in the bounded
domain.

\begin{Theorem}\label{th:plocal} Let $\Omega\subset\mathbb{R}^N$ bounded. For every $f\in L^\infty(\Omega)$ there exists
a weak solution $(u,w)$ to Problem \eqref{pp:local} with
$u(\cdot,t)\in L^\infty(\Omega)$ for every $t>0$ and $w\in
L^2([0,T];X_0^\sigma(C_\Omega))$. Moreover, the following
contractivity property holds: if $(u,w)$, $(\widetilde u,\widetilde
w)$ are the constructed weak solutions corresponding to initial data
$f,\,\widetilde f$, then
\begin{equation}
\int_{\Omega}[u(x,t)-\widetilde u(x,t)]_+\,dx\le
\int_{\Omega}[f(x)-\widetilde f(x)]_+\,dx. \label{eq:contract-p}
\end{equation}
In particular a comparison principle for constructed solutions is obtained.
\end{Theorem}

\noindent{\it Proof. } Crandall-Liggett's result only provides us in
principle with an abstract type of solution called {\sl mild
solution.} However, we know that \eqref{weak.local.bdd} and
\eqref{energy.bdd} hold, from where it is standard to show that the
mild solution is in fact weak, see for
example~\cite{dePablo-Quiros-Rodriguez-Vazquez}. We recall the main
idea: For each time $T>0$ we divide the time interval $[0,T]$ in $n$
subintervals. Letting $\varepsilon=T/n$, we construct the pair
function $(u_\varepsilon,w_\varepsilon)$ piecewise constant in each
interval $(t_{k-1},t_k]$, where $t_k=k\varepsilon$, $k=1,\cdots,n$,
as the solutions to the discretized problems
$$
\left\{
\begin{array}{ll}
L_\sigma w_{\varepsilon,k}=0\qquad &\mbox{ in } C_\Omega,\\
\varepsilon \dfrac{\partial w_{\varepsilon,k}}{\partial
y^\sigma}=w^{1/m}_{\varepsilon,k}- u_{\varepsilon,k-1}\qquad&\mbox{
on } \Omega,
\end{array}
\right.
$$
with $u_{\varepsilon,k-1}=w^{1/m}_{\varepsilon,k-1}(\cdot,0)$,
$u_{\varepsilon,0}=f$. The mild solution is obtained by
letting $\varepsilon\to0$. We still have to check that it is a weak solution.

Crandall-Ligget's Theorem gives that $u_\varepsilon$ converges in
$L^1(\Omega)$ to some function $u\in C([0,\infty): L^1(\Omega))$.
Moreover, $\|w_\varepsilon\|_\infty\le \|f\|^m_\infty$. Hence,
$w_\varepsilon$ converges in the weak-$*$ topology to some function
$w\in L^\infty(C_\Omega\times[0,\infty])$. On the other hand,
\normalcolor multiplying the equation by $w_{\varepsilon,k}$,
integrating by parts, and applying Young's inequality, we obtain
$$
\begin{array}{l}
\displaystyle\dfrac1{(m+1)}\int_\Omega
|u_\varepsilon(x,t)|^{m+1}\,dx+ \mu_\sigma\int_0^T\int_{C_\Omega}
y^{1-\sigma}|\nabla
w_\varepsilon(\overline x,t)|^2\,dxdydt\\[6pt]
\displaystyle\qquad\qquad\qquad\qquad\qquad\qquad\qquad\qquad\qquad\qquad\qquad\qquad\le
\dfrac1{(m+1)}\int_\Omega |f(x)|^{m+1}\,dx.
\end{array}
$$
Passing to the limit, the  following estimate is obtained for the
weighted norm of $|\nabla w|$,
$$
\mu_\sigma\int_0^T\int_{C_\Omega} y^{1-\sigma}|\nabla w(\overline
x,t)|^2\,dxdydt\le \dfrac1{(m+1)}\int_\Omega |f(x)|^{m+1}\,dx,
$$
and therefore $w\in L^2([0,T];X_0^\sigma(C_\Omega))$. Now, choosing
appropriate test functions, as in
\cite{dePablo-Quiros-Rodriguez-Vazquez}, it follows that we can pass
to the limit in the elliptic weak formulation to get the identity of
the parabolic weak formulation.

The contractivity obtained in each step is inherited in the limit.
\qed 

\subsection{The problem in the whole space}

\begin{Theorem}\label{th:plocal2} For every $f\in L^1(\mathbb{R}^N)\cap L^\infty(\mathbb{R}^N)$ there exists
a weak solution $(u,w)$ to Problem \eqref{pp:local}. This solution
satisfies $u(\cdot,t)\in L^1(\mathbb{R}^N)\cap
L^\infty(\mathbb{R}^N)$ for every $t>0$, and $w\in
L^2([0,T]:X^\sigma(\mathbb{R}^{N+1}_+))$. Moreover, the following
contractivity property holds: if $(u,w)$, $(\widetilde u,\widetilde
w)$ are the constructed weak solutions corresponding to initial data
$f,\,\widetilde f$, then
\begin{equation}
\int_{\mathbb{R}^N}[u(x,t)-\widetilde u(x,t)]_+\,dx\le
\int_{\mathbb{R}^N}[f(x)-\widetilde f(x)]_+\,dx.
\label{eq:contract-p2}
\end{equation}
In particular a comparison principle for constructed solutions is
obtained.
\end{Theorem}
\noindent{\it Proof. } Let us comment briefly the two
constructions of the solution.

Take to begin with as domain $\Omega=B_R$, the ball of radius $R$
centered at the origin, and let $w_R(f_R)$ be the corresponding
solution to problem \eqref{pe:local} with datum
$f_R=f\cdot\chi_{B_R}$. Passing to the limit $R\to\infty$ we obtain
a weak solution $w$ to the elliptic problem in the upper half-space
$\mathbb{R}^{N+1}_+$. Now we can follow the technique described
above using the time discretization scheme thus obtaining a weak
solution, whose trace on $\{y=0\}$ we call $U=U(f)$. This is the
mild solution produced by the Crandall-Liggett theorem and as such
is unique. The contractivity property~\eqref{eq:contract-p2} follows
from~\eqref{eq:contract-p}.

As to the second construction, we use the weak solution $u_R(f_R)$
to the parabolic problem posed in the ball $B_R$, as obtained in
Theorem~\ref{th:plocal}.  In the study of this limit we first treat
the case where $f\ge 0$ and $f_R$ approximates $f$  from below.
Then,  the family of solutions $u_R(f_R)$ is monotone in $R$ and
also $u_R(f_R)\le U(f)$ follows from simple comparison. In this way
we ensure the existence of
$$
\widetilde{u}(x,t)=\lim_{R\to\infty}u_R(x,t)
$$
It is easy to prove that $\widetilde{u}$ is also  a weak solution
with initial data $f$ and $\widetilde{u}(f)\le U(f)$. The
equivalence of the two solutions depends on the already proved
uniqueness result, see the remark after
Theorem~\ref{th:uniqueness1}.

\medskip

In the general case where $f$ changes sign, we use  comparison with
the construction for $f^+$ and $f^-$ to show that the family
$u_R(f)$ is bounded uniformly and then use compactness to pass to
the limit and obtain a weak solution.  Again we end by checking that
$\widetilde{u}= U$.~\qed

\noindent{\it Remark.} Since we have uniqueness, as a byproduct of
the limit processes of construction we get the following estimates
for weak solutions with initial data $f\in L^1(\mathbb{R}^N)\cap
L^\infty(\mathbb{R}^N)$:
$$
\|u(\cdot,t)\|_{1}\le\|f\|_{1},\qquad
\|u(\cdot,t)\|_{\infty}\le\|f\|_{\infty}.
$$

\noindent{\it Remark.} In the course of the proof we have obtained a
unique weak solution to the nonlocal nonlinear elliptic problem
\begin{equation}
\label{ell-frac} u+(-\Delta)^{\sigma/2}u^m=g\qquad \mbox{in } \Omega
\end{equation}
for every $g\in L^1(\Omega)\cap L^\infty(\Omega)$,  both for
$\Omega=\mathbb{R}^N$ and for $\Omega$ a bounded domain with
homogeneous Dirichlet condition. This weak solution satisfies  $u\in
L^1(\Omega)\cap L^\infty(\Omega)$.

\section{Existence for general data}\label{sect-L1}
\setcounter{equation}{0}

We prove here existence for data $f\in L^1(\mathbb{R}^N)$. The idea
is to approximate the initial data by a sequence of bounded
integrable functions and then pass to the limit in the approximate
problems. The key tools needed to pass to the limit are the
$L^1$-contraction property and the smoothing effect. As a
preliminary step we must show that the approximate solutions are
strong.

\subsection{Strong solutions}

We prove here that the  bounded weak  solutions constructed in the
previous section  are actually strong solutions.   We remark that
the proof does not require boundedness of the solutions, but a
control of the $L^1$ norm of the time-increment quotients. Hence the
property will be true for the general solutions constructed next by
approximation.

We start by proving something weaker, namely that $\partial_tu$ is a
Radon measure.

\begin{Proposition}
  Let $u$ be the weak solution constructed in Theorem~\ref{th:plocal2}. Then
  $\partial_t u$ is a bounded Radon measure.
\end{Proposition}

\noindent{\it Proof. } If $m=1$, a direct computation using the
representation in terms of the kernel yields
\begin{equation}
\label{eq:L1.estimate.u_t} \left\|\frac{\partial u}{\partial
t}(\cdot,t)\right\|_1\le \frac{2N}{\sigma t}\|f\|_1.
\end{equation}

If $m\ne1$, following step by step the proof for the  PME case, see
 \cite{Benilan-Crandall2} or  \cite{Vazquez},  we get
that the time-increment quotients of $u$ are bounded in
$L^1(\mathbb{R}^N)$,
\begin{equation}
\int_{\mathbb{R}^N}\dfrac1h\Big|u(x,t+h)-u(x,t)\Big|\,dx\le\frac{2}{|m-1|t}\,\|f\|_{1}+o(1)
\label{utmeasure}
\end{equation}
as $h\to0$. Hence, the limit $\partial_tu$ must be a Radon
measure.~\qed

The next step is to show that
the time derivative of a
certain power of $u$ is an $L^2_{\rm loc}$ function.

\begin{Lemma}\label{th:estimate:um+1_t} The function
$z=u^{\frac{m+1}2}$ satisfies $\partial_tz\in L^2_{\rm
loc}(\mathbb{R}^N\times(0,\infty))$.
\end{Lemma}

\noindent{\it Proof. } If we could use $\partial_t w$ as test
function, we would obtain
$$
-\mu_\sigma\int_0^t\int_{\mathbb{R}^{N+1}_+} y^{1-\sigma}\langle
\nabla(\partial w/\partial t),\nabla w\rangle\,dxdydt=
\int_0^t\int_{\mathbb{R}^N} (\partial_tz)^2dxdt,
$$
from where  we would get $\partial_tz\in
L^2(\mathbb{R}^N\times[0,T])$. Though $\partial_t w$ is not
admissible as a test function, we can work with the Steklov averages
to arrive to the same result, following \cite{Benilan-Gariepy}. For
any $g\in L^1_{\rm loc}(\mathbb{R}^N)$ we define the average
$$
g^h(x,t)=\frac1h\int_t^{t+h}g(x,s)\,ds.
$$
We see that
$$
\partial_tg^h=\frac{g(x,t+h)-g(x,t)}{h}
$$
almost everywhere. Since for our solution we have $\partial_tu^h\in
L^1(\mathbb{R}^N)$, we can write the weak formulation of solution in
the form
$$
\int_0^t\int_{\mathbb{R}^N}\partial_tu^h\varphi\,dxds=-\mu_\sigma\int_0^t\int_{\mathbb{R}^{N+1}_+}y^{1-\sigma}\langle\nabla
w^h,\nabla \varphi\rangle\,dxdyds.
$$
Let now choose the test function, $\varphi=\zeta\partial_tw^h$, where
$\zeta=\zeta(t)\ge0$. Then the above identity becomes
$$
\begin{array}{rl}
\displaystyle\int_0^t\int_{\mathbb{R}^N}\zeta\partial_t
u^h\partial(u^m)^h\,dxds&\displaystyle=-\mu_\sigma\int_0^t\int_{\mathbb{R}^{N+1}_+}y^{1-\sigma}\zeta\langle\nabla
w^h,\nabla \partial_tw^h\rangle\,dxdyds
\\ [3mm]&\displaystyle=\frac12\mu_\sigma\int_0^t\int_{\mathbb{R}^{N+1}_+}y^{1-\sigma}\zeta'|\nabla
w^h|^2\,dxdyds\le C.
\end{array}
$$
We end by using the inequality $(u^m)^h\,u^h\ge c\, [(u^{\frac{m+1}2})^h]^2$,
see for instance
\cite{dePablo-Quiros-Rodriguez-Vazquez}, and passing to the limit
$h\to0$.
 \qed

We finally  prove that $u$ is an $L^1_{\rm loc}$ function, and hence
that $\partial_tu\in L^\infty((\tau,\infty):L^1(\mathbb{R}^N))$ for
all $\tau>0$. Therefore, $u$ is a strong solution.

\begin{Theorem}\label{th:estimate:u_t} If $u$ is a weak solution to Problem~\eqref{eq:main}
such that $\partial_t u$ is a Radon measure, then $u$ is a strong solution.
\end{Theorem}

\noindent{\it Proof. } The first step is to prove that $\partial_t
u\in L^1_{\rm loc}(\mathbb{R}^N\times(0,\infty))$.   This follows
immediately from Theorem~1.1 in \cite{Benilan-Gariepy}, once we know
that $\partial_t\left(u^{\frac{m+1}2}\right)\in L^2_{\rm
loc}(\mathbb{R}^N\times(0,\infty))$, see
Lemma~\ref{th:estimate:um+1_t}. Having proved  that $\partial_t u$
is a function, estimate~\eqref{utmeasure} implies
$$
\left\|\dfrac{\partial u}{\partial t}(\cdot,t)\right\|_{1}
\le\frac{2}{|m-1|t}\,\|f\|_{1},\qquad m\ne1,
$$
while we have the estimate \eqref{eq:L1.estimate.u_t} for $m=1$.~\qed

\medskip

We end this subsection with two more estimates that will be useful in the
sequel.

Multiplying the equation by $u^m$ and integrating in space and time,
(recall that $u$ is a strong solution), we obtain
\begin{equation}
\int_0^t\int_{\mathbb{R}^N} |(-\Delta)^{\sigma/4} u^m|^2\,dxds+
\frac1{m+1}\int_{\mathbb{R}^N}|u|^{m+1}(x,t)\,dx=\frac1{m+1}\int_{\mathbb{R}^N}|f|^{m+1}\,dx.
\label{m+1}
\end{equation}
Thus, we control  the norm in $L^2((0,\infty):\dot H^{\sigma/2}(\mathbb{R}^N))$
of $u^m$ in terms of the initial data.

\begin{Proposition}
In the hypotheses of Theorems~\ref{th:general} or \ref{th:general2}, if $f\in L^{m+1}(\mathbb{R}^N)$, then the solution to Problem~\eqref{eq:main} satisfies
\begin{equation}
\int_0^\infty\|u^m(\cdot,t)\|_{\dot H^{\sigma/2}}^2\,dt\le \frac
1{m+1}\|f\|_{m+1}^{m+1}.
\label{L2grad-m+1}\end{equation}
\end{Proposition}

Another easy consequence of \eqref{m+1} is that the norm $\|u(\cdot,t)\|_{m+1}$
is nonincreasing in time. In fact, this is true for all  $L^p$ norms.
\begin{Proposition}
In the hypotheses of Theorems~{\rm \ref{th:general}} or {\rm \ref{th:general2}}, any $L^p$ norm, $1\le p\le\infty$, of the solution to Problem~\eqref{eq:main} is nonincreasing in time.
\label{pro:decr}\end{Proposition}

\noindent{\it Proof. } We multiply the equation by $|u|^{p-2}u$ with
$p>1$, and integrate in $\mathbb{R}^N$. Using Strook-Varopoulos
inequality \eqref{eq:strook.varopoulos}, we get
$$
\begin{array}{rl}
\displaystyle\frac{d}{dt}\int_{\mathbb{R}^N} |u|^p(x,t)\,dx&\displaystyle= -p\int_{\mathbb{R}^N}
(-\Delta)^{\sigma/2}(|u|^{m-1}u)|u|^{p-2}u
 \\
[4mm] &\le\displaystyle
-C\int_{\mathbb{R}^N}\left|(-\Delta)^{\sigma/4}|u|^{\frac{m+p-1}{2}}\right|^2\le0.
\end{array}$$
The limit cases $p=1$ and $p=\infty$ were obtained from the elliptic
approximation. \qed

\noindent{\it Remark.}  The previous result can be easily
generalized substituting the power $|u|^p$ by any nonnegative convex function
$\psi(u)$, using \eqref{eq:strook.varopoulos2}. Thus we obtain
$$
\frac{d}{dt}\int_{\mathbb{R}^N}
\psi(u)(x,t)\,dx
\le-\int_{\mathbb{R}^N}\left|(-\Delta)^{\sigma/4}\Psi(u)\right|^2\le0,
$$
where $\Psi(u)=\int_0^{|u|}\sqrt{ms^{m-1}\psi''(s)}\,ds$.

\subsection{Smoothing effect}\label{sect-smooth}

As a first step we obtain a bound of the $L^\infty$ norm in terms of
the $L^p$ norm of the initial datum for every $p>1$, with the
additional condition $p>p_*(m)=(1-m)N/\sigma$ if $0<m<m_*$. The
important observation, that will be used in the next subsection when
passing to the limit for general data, is that the estimates do no
depend qualitatively on the $L^\infty$ norm of the initial value.

\begin{Theorem}\label{th:smoothing2} Let $0<\sigma<2$, $m>0$, and take $p>\max\{1,(1-m)N/\sigma\}$.
Then for every $f\in L^1({\mathbb{R}^N})\cap  L^\infty({\mathbb{R}^N})$, the
solution to Problem \eqref{eq:main} satisfies
\begin{equation}
\sup_{x\in{\mathbb{R}^N}}|u(x,t)|\le C\,t^{-\gamma_p
}\|f\|_{p}^{\delta_p} \label{eq:L-inf-p2}\end{equation} with
$\gamma_p=(m-1+{\sigma}p /N)^{-1}$ and $\delta_p=\sigma p\gamma_p/N$,  the
constant $C$ depending on $m,\,p,\,\sigma$ and $N$.
\end{Theorem}

\noindent{\it Proof. } We use a classical parabolic Moser iterative
technique.

Let $t>0$ be fixed, and consider the sequence of times $t_k=(1-2^{-k})t$. We
multiply the equation in \eqref{eq:main} by $|u|^{p_k-2}u$, $p_k\ge p_0>1$, and
integrate in ${\mathbb{R}^N}\times [t_{k},t_{k+1}]$. As in the proof of
Proposition~\ref{pro:decr}, using \eqref{eq:strook.varopoulos} and the decay of
the $L^p$ norms we get
$$
\begin{array}{rl}
\displaystyle \|u(\cdot,t_k)\|_{p_k}^{p_k}&\displaystyle\ge
\frac{4mp_k(p_k-1)}{(p_k+m-1)^2}\int_{t_{k}}^{t_{k+1}}\|(-\Delta)^{\sigma/4}
|u|^{\frac{p_k+m-1}{2}}(\cdot,\tau)\|_2^2\,d\tau
\\ [4mm] &\displaystyle\ge
\frac{1}{d_k\|u(\cdot,t_k)\|_{p_k}^{p_k}}\int_{t_{k}}^{t_{k+1}}\|u(\cdot,\tau)\|_{p_k}^{p_k}
\|(-\Delta)^{\sigma/4}
|u|^{\frac{p_k+m-1}{2}}(\cdot,\tau)\|_2^2\,d\tau.
\end{array}
$$
The constant $d_k$ depends on $p_0$ (as well as on $m$ and $N$, but
not on $\sigma$). We now use the Nash-Gagliardo-Nirenberg type
inequality \eqref{eq:gagliardo.nirenberg.type.inequality} with
$r=p_k+m-1$ an use again the decay of the $L^p$ norms to obtain
$$
\begin{array}{rcl}
\displaystyle
\int_{t_{k}}^{t_{k+1}}\|u(\cdot,\tau)\|_{p_k}^{p_k}\|(-\Delta)^{\sigma/4}
|u|^{\frac{p_k+m-1}{2}}(\cdot,\tau)\|_2^2\,d\tau&\ge&\displaystyle
C\int_{t_{k}}^{t_{k+1}}\|u(\cdot,\tau)\|_{\frac{N(2p_k+m-1)}{2N-\sigma}}^{2p_k+m-1}\,d\tau
\\
[4mm] \displaystyle&\ge&\displaystyle
C2^{-k}t\|u(\cdot,t_{k+1})\|_{\frac{N(2p_k+m-1)}{2N-\sigma}}^{2p_k+m-1}.
\end{array}
$$

Summarizing, we have
$$
\|u(\cdot,t_{k+1})\|_{p_{k+1}}\le
(2^kd'_kt^{-1})^{\frac
s{2p_{k+1}}}\|u(\cdot,t_{k})\|_{p_k}^{\frac{sp_k}{p_{k+1}}},
$$
where $p_{k+1}=s(p_{k}+\frac{(m-1)}2)$, $s=\frac{2N}{2N-\sigma}>1$.

First of all we observe that taking as starting exponent $p_0=p>(1-m)N/\sigma$
(and $p>1$) it is easy to obtain the value of the sequence of exponents,
$$ p_k=A(s^k-1)+p,\qquad
A=p-\frac{(1-m)N}\sigma>0.
$$
In particular we get $p_{k+1}>p_{k}$, with $\lim_{k\to\infty}p_k=\infty$.
Observe also that $\min\{1,m\}\le \frac{p_k}{p_k+m-1}\le\max\{1,m\}$. This
implies that the coefficient in the above estimate can be bounded by
$c^{\frac{k}{p_{k+1}}}$, for some $c=c(m,p,N,\sigma)$. Now, if we denote
$U_k=\|u(\cdot,t_k)\|_{p_k}$, we have
$$
U_{k+1}\le
c^{\frac{k}{p_{k+1}}}t^{-\frac{s}{2p_{k+1}}}U_{k}^{\frac{sp_{k}}{p_{k+1}}}.
$$
This implies $ U_k\le c^{\alpha_k}t^{-\beta_k}U_{0}^{\delta_k}$,
with the exponents
$$
\alpha_k=\frac1{p_k}\sum_{j=1}^{k-1}
(k-j)s^j\to \frac{N(N-\sigma)}{\sigma^2A},\quad
\beta_k=\frac1{2p_k}\sum_{j=1}^k s^j\to\frac N{A\sigma},\quad
\delta_k=\frac{s^kp}{p_k}\to\frac{p}A\,.
$$

We conclude that
$$
\|u(\cdot,t)\|_\infty=\lim_{k\to\infty}U_k\le Ct^{-\frac N{A\sigma}}U_0^{\frac pA}
=Ct^{-\frac{N}{(m-1)N+p\sigma}}\|f\|_{p}^{\frac{p\sigma}{(m-1)N+p\sigma}}.
$$
 \qed

\noindent{\it Remark.} When $\sigma<N$, which is always the case if
$N\ge 2$, we may use the Hardy-Littlewood-Sobolev's inequality~\eqref{HLS}
instead of the Nash-Gagliardo-Nirenberg type inequality
\eqref{eq:gagliardo.nirenberg.type.inequality} to arrive at the same
result.

The constant in the previous calculations blows up both as $p\to p_*(m)$,
$0<m\le m_*$, or as $p\to1^+$ in the case $m>m_*$. Nevertheless, in this last
case an iterative interpolation argument allows to obtain the desired
$L^1$-$L^\infty$ smoothing effect.

\begin{Corollary}
  Let $0<\sigma<2$, $m>m_*$. Then for
every $f\in L^1({\mathbb{R}^N})\cap  L^\infty({\mathbb{R}^N})$, the solution to
Problem
\eqref{pp:local} satisfies
\begin{equation}
\sup_{x\in{\mathbb{R}^N}}|u(x,t)|\le C\,t^{-\gamma
}\|f\|_{1}^{\delta} \label{eq:L-inf-pp}\end{equation} with
$\gamma=\gamma_1=(m-1+{\sigma} /N)^{-1}$ and $\delta=\delta_1=\sigma\gamma/N$, the
constant $C$ depending on $m$,  $N$ and $\sigma$.
\end{Corollary}

\noindent{\it Proof. }
Putting $\tau_k=2^{-k}t$, estimate
\eqref{eq:L-inf-p2} with $p=2$ for instance (for which it is valid if $m>m_*$),
applied in the interval $[\tau_1,\tau_0]$ gives
$$
\|u(\cdot,t)\|_\infty\le c\,(t/2)^{-\gamma_2
}\|u(\cdot,\tau_{1})\|_2^{\frac{2\sigma\gamma_2}N}\le c\,(t/2)^{-\gamma_2
}\|u(\cdot,\tau_{1})\|_1^{\frac{\sigma\gamma_2}N}
\|u(\cdot,\tau_{1})\|_\infty^{\frac{\sigma\gamma_2}N}.
$$
We now apply the same estimate in the interval $[\tau_{2},\tau_{1}]$, thus
getting
$$
\|u(\cdot,t)\|_\infty\le C\,(t/2)^{-\gamma_2
}\|u(\cdot,\tau_{1})\|_1^{\frac{\sigma\gamma_2}N}
\left(C\,(t/4)^{-\gamma_2
}\|u(\cdot,\tau_{2})\|_2^{\frac{2\sigma\gamma_2}N}\right)^{\frac{\sigma\gamma_2}N}.
$$
Iterating this calculation in $[\tau_k,\tau_{k-1}]$, using
Proposition~\ref{pro:decr}, we obtain
$$
\|u(\cdot,t)\|_\infty\le C^{a_k}2^{b_k}
t^{-d_k}\|u(\cdot,0)\|_1^{e_k}\|u(\cdot,\tau_k)\|_2^{f_k}.
$$
Using the fact that  $m>m_*$ implies $\frac{\gamma_2\sigma}{
N}=\frac{\sigma}{(m-1)N+2\sigma}<1$, we see that the exponents satisfy, in the
limit $k\to\infty$,
$$
\begin{array}{l}
\displaystyle
a_k=\sum_{j=0}^{k-1}\Big(\frac{\gamma_2\sigma}N\Big)^j\to\frac{(m-1)N+2\sigma}{(m-1)N+\sigma}=
\frac{\gamma_1\sigma}N+1\,, \\
[4mm]
\displaystyle
b_k=\sum_{j=0}^{k-1}\gamma_2 (j+1)\Big(\frac{\gamma_2\sigma}N\Big)^j\to
\frac{(m-1)N+2\sigma}{((m-1)N+\sigma)^2}\,, \\
[4mm]
\displaystyle
d_k=\gamma_2 a_k\to\frac{N}{(m-1)N+\sigma}=\gamma_1\,,
\\ [4mm]
\displaystyle
e_k=a_k-1\to \frac{\gamma_1\sigma}N\,,
\\ [4mm]
\displaystyle
f_k=2\Big(\frac{\gamma_2\sigma}N\Big)^k\to 0\,.
\end{array}
$$
\qed

\noindent{\it Remark.}  Since the $L^p$ norm is nonincreasing, an $L^r$ decay is
obtained again by interpolation, for any $r\ge p\ge1$,
$$
\|u\|_r\le ct^{-\frac{(r-p)N}{r((m-1)N+p\sigma)}}\,\|f\|_p^{\frac{p((m-1)N+r\sigma)}{r((m-1)N+p\sigma)}}.
$$

\subsection{Passing to the limit. Existence of strong solutions}
We first prove the existence of a strong solution for the case of initial data with improved regularity.
\begin{Theorem}
\label{finalexistence} Let $\sigma\in (0,2)$  and  $m>0$. Then for
every  $f\in L^1(\mathbb{R}^N)$ if $m>m_*$ or $f\in
L^1(\mathbb{R}^N)\cap L^p(\mathbb{R}^N)$  with
$p>p_*(m)=(1-m)N/\sigma$ if $m\le m_*$,
 there exists a  strong solution
 to Problem~\eqref{eq:main}.
\end{Theorem}

\noindent{\it Proof.} Let $\{f_k\}\subset L^1({\mathbb{R}^N})\cap
L^\infty({\mathbb{R}^N})$ be a sequence of functions converging to
$f$ in $L^1$,  and let $\{u_k\}$ be the sequence of the
corresponding solutions. Thanks to the $L^1$-contraction property,
we know that $u_k(\cdot,t)\to u(\cdot,t)$ in $L^1({\mathbb{R}^N})$
for all $t>0$ for some function $u$. Moreover, nonlinear Semigroup
Theory guarantees that $u_k\to u$ in
$C([0,\infty):L^1(\mathbb{R}^N))$ \cite{Crandall}, \cite{Evans}.

Consider a fixed time $\tau>0$. Using the smoothing effect and the
estimate~\eqref{L2grad-m+1}, we have $u_k\in
L^\infty({\mathbb{R}^N}\times[\tau,\infty))$ and $u^m_k\in
L^2((\tau,\infty):\dot H^{\sigma/2}({\mathbb{R}^N}))$, both
uniformly in $k$. Thus the limit $u$ is a weak solution to
Problem~\eqref{eq:main} for every $t\ge\tau$. We now want to go down
to $\tau=0$. This follows from the $L^1$-contraction and the
$L^1$-continuity. In fact
$$
\begin{array}{rl}
\displaystyle\int_{\mathbb{R}^N}|u(x,t)-f(x)|\,dx&\displaystyle\le\int_{\mathbb{R}^N}|u(x,t)-u_k(x,t)|\,dx+
\int_{\mathbb{R}^N}|u_k(x,t)-f_k(x)|\,dx \\ [4mm]
&\displaystyle+\int_{\mathbb{R}^N}|f_k(x)-f(x)|\,dx.
\end{array}
$$
The fact that $u$ is a strong solution is now an immediate consequence of Theorem~\ref{th:estimate:u_t}.~\qed

As a byproduct, and using the uniqueness results, Theorems~\ref{th:uniqueness1}
and~\ref{th:uniqueness2}, we obtain the $L^1$-ordered contraction for the limit
solution: given $u$ and $\widetilde u$ two weak solutions to Problem
\eqref{eq:main}, then for every $0\le t_1<t_2$ we have
\begin{equation}
\int_{{\mathbb{R}^N}}[u(x,t_2)-\widetilde u(x,t_2)]_+\,dx\le
\int_{{\mathbb{R}^N}}[u(x,t_1)-\widetilde u(x,t_1)]_+\,dx.
\label{eq:contract-l1}\end{equation}

We now consider the case of general data, $f\in L^1(\mathbb{R}^N)$.
We only need to look at the case $0<m\le m_*$, in view of the
previous theorem. We prove here that the unique mild solution is in fact a very
weak solution.
\begin{Theorem}
\label{th:mild.very.weak}
  Let $\sigma\in (0,2)$  and  $m\le m_*$. Then for every  $f\in L^1(\mathbb{R}^N)$ there is a unique mild solution to Problem~\eqref{eq:main} which is moreover a very weak solution.
\end{Theorem}

\noindent{\it Proof.} As in the proof of the previous theorem, we
approximate the initial data by a sequence $\{f_k\}$ of bounded, integrable
initial data. The corresponding strong solutions $u_k$ converge in
$C([0,\infty):L^1(\mathbb{R}^N))$ to a certain function $u$ which, being the
limit of mild solutions,  is also a mild solution  according to the theory,
\cite{Crandall-Liggett}.
 Moreover, integrating by parts in space, we see that the approximate solutions $u_k$ are very weak solutions, namely
\begin{equation}
\int_0^\infty\int_{\mathbb{R}^N}u_k\,\frac{\partial \varphi}{\partial t}-\int_0^\infty\int_{\mathbb{R}^N}u_k^m(-\Delta)^{\sigma/2}\varphi=0
\end{equation}
for every smooth and compactly supported test function $\varphi$. On the other
hand, $\|u_k^m\|_{1/m}=\|u_k\|_{1}^m \le \|f\|_1^m$. Hence
$u_k^m\rightharpoonup u^m$ in $L^{1/m}(\mathbb{R}^N)$. Since
$(-\Delta)^{\sigma/2} \varphi$ belongs to the dual space
$L^{1/(1-m)}(\mathbb{R}^N)$, we conclude that $u$ is a very weak solution.~\qed

The passage to the limit in the case where the spatial domain is bounded is similar.


\section{Further qualitative properties of the solutions}\label{sect-props}
\setcounter{equation}{0}

We prove in this section some important properties that our solutions have.
Throughout this section $u$ is the strong solution to Problem~\eqref{eq:main}
corresponding to an initial value $f$ satisfying the hypotheses of
Theorem~\ref{finalexistence}.

\subsection{Positivity and regularity}

We start with this observation: if $u\ge0$ is a classical solution
and $u(x_0,t)=0$ for some $x_0$ and $t$, then formula
\eqref{def-riesz} gives $(-\Delta)^{\sigma/2}u^m(x_0,t)<0$, unless
$u(\cdot,t)\equiv0$, and hence $\partial_t u(x_0,t)>0$. Therefore,
we expect solutions with nonnegative data to become positive
immediately, and to stay positive unless they vanish. However,
solutions are not known to be classical. Hence, we  will use a
different argument, which involves the extension
Problem~\eqref{pp:local}.

The first ingredient in our proof  is a control of the decay of the solutions,
that has an independent interest. This control is based in an argument of
Alexandrov's type, cf.~\cite{ArCaff8}, \cite{Vazquez}, which is a bit delicate in this
case, since the function $u$ is not yet known to be continuous.

\begin{Proposition}\label{alex-nocont}
Assume $f$ has compact support. For any bounded measurable set $M$
with $|M|>0$ there is a large enough radius $R_*$, depending only on
the support of the initial data and on $M$, such that
$$
u(x,t)\le \sup_{z\in M}u(z,t)\qquad\text{ for a.e. } |x|>R_*,\; t>0.
$$
\end{Proposition}
\noindent{\it Proof. } Let $\mbox{supp}(f)\subset B(0,R)$, $M\subset
B(0,R')$. Thanks to Lebesgue's density Theorem, we know that there
is a point $x_1\in M$ such that for all $\delta>0$ there exists a
radius $r_\delta$ such that for all cubes $Q\subset B(x_1,r_\delta)$
we have
$$
\frac{|M\cap Q|}{|Q|}>1-\delta.
$$

Let us now consider the cube $Q^*=Q(x^*,2)$, with $x^*$ far away
from the origin to be chosen later. We can cover $Q^*$ (except a
subset of zero measure) with a finite number of disjoint cubes,
small enough such that they are reflections of cubes centered at
$x_1$ and contained in $B(x_1,r_\delta)$. Hence, an argument of
Alexandrov's type (which can be done if $|x|^*$ is large enough)
shows that
$$
\frac{|\{u(\cdot,t)<\sup_{x\in M}u(x,t)\}\cap Q^*|}{|Q^*|}>1-\delta.
$$
More precisely, let $\widetilde Q$ be any of the cubes covering
$Q^*$, and let $\widetilde x$ be its center. We reflect around the
hyperplane in $\Omega$, $\pi\equiv(\widetilde
x-x_1)\cdot(x-(x_1+\widetilde x)/2)=0$, getting a cube $Q'\subset
B(x_1,r_\delta)$. It is clear that if $|\widetilde x|$ is large
(depending on $R$ and $R'$), then the hyperplane $\pi$ divides the
half-space $\overline\Omega$ in two parts, $\overline\Omega=H_1\cup
H_2$ with $B(0,R)\times[0,\infty)\subset H_1$, $\widetilde
Q\times[0,\infty)\subset H_2$. In this way, by the comparison
principle, we obtain that the function
$z(x,y,t)=w(x,y,t)-w(\widetilde x+x_1-x,y,t)$ satisfies
$z(x,y,t)\ge0$ almost everywhere in $H_1$, $t>0$.

Therefore,
$$
|\{u(\cdot,t)<\sup_{x\in M}u(x,t)\}\cap \widetilde
Q|\ge|\{u(\cdot,t)<\sup_{x\in M}u(x,t)\}\cap
Q'|>(1-\delta)|Q'|=(1-\delta)|\widetilde Q|.
$$
Summing up,
$$
|\{u(\cdot,t)<\sup_{x\in M}u(x,t)\}\cap Q^*|>(1-\delta)|Q^*|.
$$

We end by letting $\delta\to0$ to get that $u(\cdot,t)<\sup_{x\in
M}u(x,t)$ a.e.~in $Q^*$. \qed

The second ingredient in the proof of positivity is  the following
estimate on the time derivative for nonnegative solutions, arising
from the homogeneity of the parabolic operator, see for instance
\cite{Benilan-Crandall},
\begin{equation}
(m-1)t\frac{\partial u}{\partial t}+ u \ge 0.
\label{wt+w>}\end{equation}

If $m=1$ this formula is empty. However, in that case
 we have the representation formula \eqref{lineal}, in
terms of the fundamental solution $K_\sigma$, from which  it is easy
to derive the estimate
$$
\sigma t\partial_t u+ Nu \ge 0.
$$
It is sharp: equality holds for the fundamental solution at the
origin.

\begin{Theorem}\label{th:positivity}
Let $f\ge0$. Then for every $t>0$: either $u(\cdot,t)\equiv0$ or $
\inf_{K}u(x,t)\ge c_K>0$ for every compact $K\subset{\mathbb{R}^N}$.
\end{Theorem}

\noindent{\it Proof. } By comparison, we only need to consider
compactly supported initial data. Let then the support of $f$ be contained in
the ball $B(0,R)$. Let $K\subset{\mathbb{R}^N}$ be any compact set, and assume
that there exists some $t_0>0$ such that $
\inf_{K}u(x,t_0)=0$. Then there is a measurable set
$M_\varepsilon\subseteq K$, $|M_\varepsilon|>0$, such that
$u(\cdot,t_0)<\varepsilon$ in $M_\varepsilon$. By
Proposition~\ref{alex-nocont}, there exists some ball $B(x^*,1)$
where $u(\cdot,t_0)<\varepsilon$ almost everywhere. Letting
$\varepsilon\to0$ we get  $u(\cdot,t_0)=0$ almost everywhere in
$B(x^*,1)$. Then we use \eqref{wt+w>} to prove that $u$ is zero in
that ball for an interval of times $I$. In fact, if $m>1$ we obtain
from \eqref{wt+w>}
$$
u(\cdot,t_2)\ge u(\cdot,t_1)\Big(\frac{t_1}{t_2}\Big)^{1/(m-1)},
$$
a.e.~in ${\mathbb{R}^N}$, $t_2\ge t_1>0$. Therefore we deduce $u(\cdot,t)=0$
a.e.~in $B(x^*,1)$ for every $0<t\le t_0$, that is, we may take $I=[0,t_0]$. In
the case $m<1$ formula \eqref{wt+w>} gives the same property in
$I=[t_0,\infty)$. The case $m=1$ follows from
\eqref{lineal}.

We now observe that the integral definition of solution implies that
for every test function $\varphi$ that vanishes on $\partial
B(x^*,1)\times(0,\infty)$, we have
$$
\begin{array}{rl}
0=&\displaystyle\int_I\int_{B(x^*,1)}
u\frac{\partial\varphi}{\partial t}\,dxdt=
\mu_\sigma\int_I\int_0^\infty\int_{B(x^*,1)}
y^{1-\sigma}\langle\nabla
w,\nabla\varphi\rangle\,dxdydt \\
[3mm] =&\displaystyle-\int_I\int_{B(x^*,1)}\frac{\partial
w}{\partial y^\sigma}\,\varphi\,dxdt.
\end{array} $$
This gives $\partial_{y^\sigma} w (\cdot,0,t)=0$ a.e.~in $x\in
B(x^*,1)$ for all $t\in I$.  Now, for each fixed $t\in I$ we extend
$w$ in a even way in the $y$ variable to obtain a solution of the
elliptic equation $L_\sigma w=0$ in a ball
$B((x^*,0),1)\subset\mathbb{R}^{N+1}$. Since
  $A(x,y)=|y|^{1-\sigma}$ is an
$A_2$-weight, we can apply  a half Harnack's inequality (Theorem
2.3.1 in \cite{Fabes-Kenig-Serapioni}), to obtain
$$
\inf_{B(x^*,1/2)}u^m(x,t)\ge\inf_{B((x^*,0),1/2)}w(\overline x,t)\ge
c\|w(\cdot,t)\|_{L^2(B((x^*,0),1))}.
$$
If $u(\cdot,t)\not\equiv0$ in ${\mathbb{R}^N}$, then $w(\cdot,t)>0$ in
$\Omega$, and this results in a contradiction. \qed

As a consequence of positivity, we can extend to the whole range
$m>0$ the continuity of solutions that was obtained in
\cite{Athanasopoulos-Caffarelli} for $m\ge 1$.

\begin{Theorem}\label{th:regularity2}
Assume $u\ge0$ and $u(\cdot,T)\not\equiv0$. Then $u\in
C^\alpha({\mathbb{R}^N}\times(0,T))$ for some $0<\alpha<1$.
\end{Theorem}

\noindent{\it Proof. } The above-mentioned regularity result of
\cite{Athanasopoulos-Caffarelli} applies for bounded solutions to the equation
$$
\dfrac{\partial \beta(v)}{\partial t} + (-\Delta)^{\sigma/2} v=0
$$
in some ball $B\subset{\mathbb{R}^N}$ and $t>0$, with a nondegeneracy condition
on the constitutive monotone function $\beta$. This condition is fulfilled once
we know that in any given ball the solution is essentially bounded below away
from zero. On the other hand the solution is bounded for every positive time,
thanks to the smoothing effect.~\qed

\subsection{Conservation of mass}

\begin{Theorem}\label{th:mass}
Let $m\ge m_*$. Then for every $t>0$ we have
$$
\int_{{\mathbb{R}^N}}u(x,t)\,dx= \int_{{\mathbb{R}^N}}f(x)\,dx.
$$
\end{Theorem}

\noindent{\it Proof. } Thanks to the $L^1$-contraction property, it
is enough to consider the case of bounded initial data. The proof will follow
different arguments in the cases $m>m_*$ and $m=m_*$, and even in this latter
case we have to distinguish between dimensions $N\ge2$ and $N=1$.

\noindent{\sc Case $m>m_*$.} We can adapt the technique that was used to prove
the property in the  case $\sigma=2$ to deal with the nonlocal operator. It
works as follows: we take a nonnegative non-increasing cut-off function
$\psi(s)$ such that $\psi(s)=1$ for $0\le s\le1$, $\psi(s)=0$ for $s\ge2$, and
define $\varphi_R(x)=\psi(|x|/R)$. Multiplying the equation by $\varphi_R$  and
integrating by parts, we obtain, for every $t>0$,
\begin{equation}
\displaystyle\frac{d}{dt}\int_{\mathbb{R}^N}
u\varphi_R=-\int_{{\mathbb{R}^N}}u^m\,(-\Delta)^{\sigma/2}\varphi_R.
\label{eq:mass-cons}\end{equation}

The radial cut-off function $\varphi_R$ has the scaling property
\begin{equation}
\label{varphi1}
(-\Delta)^{\sigma/2}\varphi_R(x)=R^{-\sigma}(-\Delta)^{\sigma/2}\varphi_1(x/R).
\end{equation}
In addition,
 $(-\Delta)^{\sigma/2}\varphi_1\in L^1(\mathbb{R}^N)\cap
L^\infty(\mathbb{R}^N)$. Both properties are straightforward using
representation~\eqref{def-riesz}.

Then, if we apply H\"older's inequality with $p=\max\{1,1/m\}$ to the right-hand side of
\eqref{eq:mass-cons}, and use the above property, we get
$$
\begin{array}{rl}
\displaystyle\left|\frac{d}{dt}\int_{\mathbb{R}^N} u\varphi_R\right|&\displaystyle\le
\|f\|_{\infty}^{m-1/p}
\|f\|_{1}^{1/p}\|(-\Delta)^{\sigma/2}\varphi_R\|_{p/(p-1)} \\
&\displaystyle\le R^{-\sigma+N(p-1)/p}\|f\|_{\infty}^{m-1/p}
\|f\|_{1}^{1/p}\,.
\end{array}
$$
We conclude since the exponent of $R$ is negative precisely for $m>m_*$, and
thus
$$
\frac{d}{dt}\int_{\mathbb{R}^N}
u=\lim_{R\to\infty}\frac{d}{dt}\int_{\mathbb{R}^N} u\varphi_R=0.
$$

\medskip

\noindent{\sc Case $m=m_*$, $N\ge2$.}
This case is much more difficult. The idea is to study separately the behaviour of
the mass in a bounded set and close to the infinity, and also to decompose the
fractional Laplacian into two operators. First, for every given $\delta>0$, we
put $u=u_1+u_2$, where
$$
u_1=u\cdot\chi_{\{|x|<R_0\}},\quad \int_{\mathbb{R}^N} |u_2|<\delta.
$$
Observe that $u^m=u_1^m+u_2^m$. Now express these functions in the following
form
$$
u_1^m=(-\Delta)^{\gamma/2}z,\quad u_2^m=\varepsilon w+(-\Delta)^{\gamma/2}w,
$$
where $\gamma=2-\sigma>0$, $\varepsilon>0$. Then our equation becomes
$$
\frac{\partial u}{\partial t}=\Delta z+\Delta
w-\varepsilon(-\Delta)^{\sigma/2}w.
$$
We have introduced the $\varepsilon$-regularization in the definition of $w$
since, when applying inequality \eqref{HLS}, it gives no information in the
critical case $m=m_*$ if $\varepsilon=0$.

As before, multiplying by the test function $\varphi_R$, and integrating by
parts, we have
$$
\begin{array}{rl}
\displaystyle\left|\frac{d}{dt}\int_{\mathbb{R}^N} u\varphi_R\right|
&\displaystyle\le
\int_{\mathbb{R}^N}|z\Delta\varphi_R|+\int_{\mathbb{R}^N}|w\Delta\varphi_R|+
\varepsilon\int_{\mathbb{R}^N}|(-\Delta)^{\gamma/2}w\,(-\Delta)^{\sigma-1}\varphi_R| \\
[3mm] &\displaystyle= I_1+I_2+I_3\,.
\end{array}
$$
We estimate each integral using H\"older's inequality, the properties of
$\varphi_R$ and some estimates on $z$ and $w$.

\noindent{\it Estimate of $I_1$}. We have $u_1^m\in L^r({\mathbb{R}^N})$ for
every $1\le r\le\infty$, with
$$
\|u_1^m\|_r\le \|u_1\|_1^{mr}|\{|x|<R_0\}|^{\frac{(1-mr)}{r}}\le CR_0^{\frac{N(1-mr)}{r}}\,.
$$
Using  inequality \eqref{HLS} we have that $z=(-\Delta)^{-\gamma/2}(u_1^m)\in
L^q({\mathbb{R}^N})$, $q=\dfrac{Nr}{N-\gamma r}$, for every $1<r<N/\gamma$,
with $\|z\|_q\le c\|u_1^m\|_r$. Then for such values of $r$, and using
H\"older's inequality and \eqref{varphi1}, we have
$$
|I_1|\le \|z\|_q\|\Delta \varphi_R\|_{q/(q-1)}\le c R_0^{\frac{N(1-mr)}{r}}
R^{-2+N(q-1)/q}=c(R_0/R)^{\frac{N(1-mr)}{r}}.
$$
The exponent is positive if we take $1<r<1/m$.

\medskip

\noindent{\it Estimate of $I_2$}. We have here $u_2^m\in
L^r({\mathbb{R}^N})$ for every $1/m\le r\le\infty$, with
$$
\|u_2^m\|_{1/m}= \|u_2\|_1^m\le c\delta^m.
$$
If we multiply now the equation satisfied by $w$ by $w^{1/m-1}$, and integrate
in ${\mathbb{R}^N}$, we get
$$
\varepsilon\int_{\mathbb{R}^N} w^{1/m}+\int_{\mathbb{R}^N} w^{1/m-1}(-\Delta)^{\gamma/2}w=\int_{\mathbb{R}^N}
u_2^mw^{1/m-1}.
$$
The second term is nonnegative by \eqref{eq:strook.varopoulos}. Thus, by
H\"older's inequality we get
$$
\varepsilon\|w\|_{1/m}\le \|u_2^m\|_{1/m}\|w^{1/m}\|_{1/(1-m)}^{1-m},
$$
i.e., $\varepsilon\|w\|_{1/m}\le\|u_2^m\|_{1/m}\le c\delta^m$.  This implies
$$
|I_2|\le\|w\|_{1/m}\|\Delta \varphi_R\|_{1/(1-m)}\le
c\delta^m\varepsilon^{-1}R^{\sigma-2}.
$$

\noindent{\it Estimate of $I_3$}. Since from the previous
calculations we have $\|(-\Delta)^{\gamma/2}w\|_{1/m}\le2\|u_2^m\|_{1/m}\le
c\delta^m$, we get
$$
|I_3|\le
\varepsilon\|(-\Delta)^{\gamma/2}w\|_{1/m}\|(-\Delta)^{\sigma-1}\varphi_R\|_{1/(1-m)}
\le c\delta^m\varepsilon R^{2-\sigma}.
$$

\medskip

Summing up, we have obtained
$$
\left|\frac{d}{dt}\int_{\mathbb{R}^N} u\varphi_R\right|\le
c(R_0/R)^{\frac{N(1-mr)}{r}}+c\delta^m\varepsilon^{-1}R^{\sigma-2}+c\delta^m\varepsilon
R^{2-\sigma}.
$$
We now choose $\varepsilon=R^{\sigma-2}$, and make first $R\to\infty$ and then
$\delta\to0$ to conclude.

\medskip

\noindent{\sc Case $m=m_*$, $N=1$.} Observe that $m_*>0$  implies $\sigma<1$.
Here we consider the same functions $z$ and $w$ as before, but with
$\gamma=1-\sigma>0$. The equation becomes in this case
$$
\frac{\partial u}{\partial t}=-(-\Delta)^{1/2}z-(-\Delta)^{1/2}
w-\varepsilon(-\Delta)^{\sigma/2}w.
$$
From here on, the calculations are exactly the same as the ones for
the case $m=m_*$, $N\ge2$.~\qed

\subsection{Extinction}

The condition $m\ge m_*$ to have mass conservation is not technical, as shown
by  the next result on extinction in finite
 time, which extends the result by  B\'{e}nilan and Crandall
 for the standard differential case $\sigma=2$, see \cite{Benilan-Crandall}.

\begin{Theorem}
  Let $0<\sigma<\min\{2,N\}$ and $0<m<(N-\sigma)/N$. Then,
 if $f\in L^{(1-m)N/\sigma}({\mathbb{R}^N})$,  there is a finite time $T>0$ such that  $u(x,T)=0$ a.e.~in ${\mathbb{R}^N}$.
\label{pro:extinction}
\end{Theorem}

\noindent{\it Proof. } From the proof of Proposition~\ref{pro:decr},
we have
$$
\frac{d}{dt}\int_{\mathbb{R}^N} |u|^{p}\le
-C\int_{{\mathbb{R}^N}}\left|(-\Delta)^{\sigma/4}|u|^{\frac{m+p-1}{2}}\right|^2.
$$
Using  inequality \eqref{HLS}
we obtain
$$
\frac{d}{dt}\int_{\mathbb{R}^N} |u|^p\,dx+C\Big(\int_{\mathbb{R}^N}
|u|^{\frac{(p+m-1)N}{N-\sigma}}\,dx\Big)^{\frac{N-\sigma}N}\le 0.
$$
If we now choose $p=(1-m)N/\sigma$ (this is where the restriction on
$m$ comes, since $p$ has to be bigger than one), we get that the
function $J(t)=\|u(\cdot,t)\|_p$ satisfies the differential
 inequality
$$
J'(t)+CJ^{\frac{N-\sigma}N}(t)\le0.
$$
This implies extinction in finite time  provided $J(0)$ is
finite.~\qed

\noindent\emph{Remark. }  In the special case
$m=(N-\sigma)/(N+\sigma)$, $N>\sigma$,  there exists an explicit family of
solutions in separated variables with extinction in finite time,
$$
u(x,t)=b_{N,\sigma,c}(T-t)^{\frac{N+\sigma}{2N}}[c+|x-a|^2]^{-\frac{N+\sigma}2},\qquad
T,\,c>0,\,a\in\mathbb{R}^N.
$$
The spatial part satisfies the elliptic fractional equation
$(-\Delta)^{\sigma/2}\varphi^m=\varphi$, see for instance
\cite{Chen-Li-Ou}.

\subsection{The problem in a bounded domain}
It is easy to see, following the ideas of the proofs for the case
where the domain is the whole space, that the smoothing
effect is true for the solutions to the problem posed in a bounded
domain $\Omega\subset\mathbb{R}^N$ (though the obtained decay rate will  not be
optimal for the problem in a bounded domain). The decay of the solution to 0 yields the
decay of the mass. Even more, solutions become extinct in a finite
time for very $0<m<1$.

\begin{Proposition}
  Let $0<\sigma<2$ and $0<m<1$. Then,
 if $f\in L^{p}(\Omega)$ for some $p>\max\{1,(1-m)N/\sigma\}$,
 there is a finite time $T>0$ such that  $u(x,T)=0$ a.e.~in $\Omega$.
\label{pro:extinction2}
\end{Proposition}

\noindent{\it Proof. } With the same calculations used for the whole
space, using here inequality \eqref{eq:gagliardo.nirenberg.type.inequality}, we
obtain
$$
\frac{d}{dt}\int_{\Omega} |u|^p\,dx+C\Big(\int_{\Omega}
|u|^{p}\,dx\Big)^{\gamma}\le 0,
$$
where $\gamma=(p+m-1)/p\in(0,1)$.~\qed

\noindent {\sc Retention.} For $m\ge1$ it is easy to see that
nonnegative solutions do not become extinct in finite time, not only
in the case of the Cauchy Problem posed in $\mathbb{R}^N$, but also
for the problem posed in a bounded domain with zero Dirichlet data.
Indeed,  if $m>1$ nonnegative solutions  satisfy
estimate~\eqref{wt+w>}. Hence, the function $ t^{\frac1{m-1}} u$ is
nondecreasing.  This implies a retention property: if a solution is
positive at some point at some time, it will remain positive at that
point for any later times.

 The above retention property for nonnegative solutions is also
true for the case $m=1$, though a different proof is needed. Let
$\mathcal{O}\subset\Omega$ be an open set in which $u(\cdot,t_1)>0$ (recall
that $u$ is continuous for $t>0$). Let $\varphi=\varphi_{\mathcal{O},1}$ be the
normalized eigenfunction corresponding  to the first eigenvalue
$\mu=\mu_{\mathcal{O},1}$ of $(-\Delta)^{\sigma/2}$ in $\mathcal{O}$ with
homogeneous Dirichlet boundary condition. We take $\varepsilon>0$ small enough
so that $u(\cdot,t_1)\ge
\varepsilon\varphi$ a.e.~in $\mathcal{O}$. The solution to the
fractional heat equation in $\mathcal{O}$, $t\ge t_1$, with initial data
$\varepsilon\varphi$ and homogeneous Dirichlet boundary condition is given by
$$
v(x,t)=\varepsilon e^{-\mu(t-t_1)}\varphi(x).
$$
We claim that $u\ge v$ for all $x\in\mathcal{O}$, $t\ge t_1$,  from
where the retention property follows, since $\varphi>0$ in
$\mathcal{O}$.

 The claim is proved using a comparison argument in
$\mathcal{O}$. Comparison is in principle not at all obvious, since
the fractional Laplacian operator changes with the domain. However,
we have
$$
(-\Delta)^{\sigma/2}_{\mathcal{O}}u\ge
(-\Delta)^{\sigma/2}_{\Omega}u\qquad\mbox{in }
\mathcal{O}\subset\Omega.
$$
This is proved easily for any nonnegative $u$ in the right spaces using the
extensions of $u$ to the corresponding half-cylinders. Hence
$\partial_tu+(-\Delta)^{\sigma/2}_{\mathcal{O}}u^m\ge0$, from where the claim
follows.

An analogous argument can be performed for nontrivial nonpositive  solutions to
prove that they never become extinct if $m\ge1$.

\medskip

\noindent{\it Remark. } In the case $m=1$ non-extinction
is true for any nontrivial initial data in $L^1(\Omega)$. Indeed,
the solution becomes bounded immediately after $t=0$, and remains
nontrivial for a while (recall that it is continuous in
$L^1(\Omega)$). In particular, $u(\cdot,\tau)\in L^2(\Omega)$ and is
nontrivial for some small $\tau$.  Since $u$ is smooth for positive
times (in fact $C^\infty$), the solution can be expanded as a series
of eigenfunctions, $u(\cdot,t)=\sum_{k=1}^\infty a_k
e^{-\mu_k(t-\tau)}\varphi_k$. Hence
$\|u(\cdot,t)\|_{L^2(\Omega)}=\sum_{k=1}^\infty|a_k|^2e^{-2\mu_k(t-\tau)}$
for all $t\ge\tau$. Notice that this quantity is positive for all
$t\ge\tau$, since $u(\cdot,\tau)$ was nontrivial.


\section{Continuous dependence}\label{sect-cont.dep.}
\setcounter{equation}{0}

In this section we prove the continuous dependence  part of
Theorem~\ref{th:properties}.  The result we obtain below includes the limit
case $m=m_*$ for $N>2$.  The case $\sigma
\to 2$ is a bit different from the rest, and is dealt with separately.

\subsection{$0<\sigma<2$}
We consider for $N>2$ the region
$$
D=\{(m,\sigma) : 0<\sigma<2, \ m\ge(N-\sigma)/N\},
$$
while for $N=1,2$ we take
$$
D=\{(m,\sigma) : 0<\sigma<2, \ m>(N-\sigma)_+/N\}.
$$

\begin{Theorem}
 The map $S: D\times L^1({\mathbb{R}^N})\to C([0,\infty):\,L^1({\mathbb{R}^N}))$ is
continuous in all the arguments $(m,\sigma, f)$. \label{th:cont}
\end{Theorem}

This  will follow from a result of nonlinear Semigroup Theory which
states that if each of $A_n$, $n=1,2,\dots,\infty$ is an
$m$-accretive operator in a Banach space $\cal X$,
$f_n\in\overline{D(A_n)}$ and $u_n$ is the solution of
$$
\frac{du_n}{dt}+A_nu_n=0,\qquad u_n(0)=f_n,
$$
then $A_n\to A_\infty$ and $f_n\to f_\infty$ imply $u_n\to u_\infty$
in $C([0,\infty):\, {\cal X})$, where $A_n\to A$ is understood as
$$
\lim_{n\to\infty}(I+A_n)^{-1}g= (I+A_\infty)^{-1}g\quad\mbox{for all
} g\in {\cal X}.
$$
See, e.g, \cite{Crandall}, \cite{Evans} for statements and references. Hence,
the theorem will be a corollary of the convergence of
$(I+A_{m_n,\sigma_n})^{-1}$, where $A_{m,\sigma}(u)=(-\Delta)^{\sigma/2}u^{m}$.
Thanks to the contractivity in $L^1({\mathbb{R}^N})$ of the elliptic problems
under consideration,  it is enough to prove this convergence for functions $g$
which are also bounded.

\begin{Proposition}
Let $\{(m_n,\sigma_n)\}_{n=1}^\infty$, $(m_n,\sigma_n)\in D$, be such that
$m_n\to \overline m$ and $\sigma_n\to\overline \sigma$ as $n\to\infty$,
$(\overline m,\overline \sigma)\in D$. Then, for all $g\in
L^1({\mathbb{R}^N})\cap L^\infty({\mathbb{R}^N})$
$$\lim\limits_{n\to\infty}(I+A_{m_n,\sigma_n})^{-1}g= (I+A_{\overline
m,\overline\sigma})^{-1}g.
$$
\end{Proposition}
\noindent{\it Proof. } {\sc Step 1. } Let
$u_{n}=(I+A_{m_n,\sigma_n})^{-1}g$, i.e., $u_n$ is the unique
solution to the equation
\begin{equation}
u_n+(-\Delta)^{\sigma_n/2}u_n^{m_n}=g,
\label{eq-dep-cont}\end{equation}
see \eqref{ell-frac}. The $L^1$-contraction estimate~\eqref{eq:contract-e}, which is also valid for the case where $\Gamma=\mathbb{R}^N$,
implies the bounds
$$
\begin{array}{l}
\|u_{n}\|_{1}\le \|g\|_{1}\,,\\
[3mm] \|u_{n}-\tau_h u_{n}\|_{1}\le\|g-\tau_hg\|_{1}\,,
\end{array}
$$
for each $h\in{\mathbb{R}^N}$, where $(\tau_h v)(x)=v(x+h)$. This is enough,
thanks to Fr\'echet-Kolmogorov's compactness criterium, to prove that
$\{u_{n}\}$ is precompact in $L^1(K)$ for each compact set
$K\subset{\mathbb{R}^N}$.

\noindent {\sc Step 2. }   To extend compactness to the whole ${\mathbb{R}^N}$ we need to
control uniformly the tails of the solutions  at infinity.  More precisely, we
need to prove that, given $\varepsilon>0$, there exists some $R>0$ such that
$\|u_{n}\|_{L^1({\mathbb{R}^N}\setminus B_{R}(0))}<\varepsilon$. This will
follow from a computation which is very similar to the one in the proof of
Theorem~\ref{th:mass}, but now taking as test function $1-\varphi_R$ instead of
$\varphi_R$. The problem is that this convergence fails to be uniform when $m$
approaches $m_*$. In order to reach $m_*$ we use a different argument, which
unfortunately only works when $N>2$.

We first perform a reduction in order to consider only a characteristic
function as initial value,  a technique borrowed from
\cite{Benilan-Crandall}. Given $\varepsilon>0$, there exist
$M,r_0>0$ such that $h=M\chi_{\{|x|\le r_0\}}$ satisfies $
\|(g-h)_+\|_{1}<\varepsilon$. Notice that $g\le g^+\le
h+(g-h)_+$. Hence, using the $L^1$-contraction property
\eqref{eq:contract-e}, we have, denoting by $u_{h}$ the solution
corresponding to the initial data $h$,
$$
\begin{array}{rl}
\displaystyle\int_{|x|>R}(u_g)^+&\displaystyle\le\int_{|x|>R}u_{g^+}\le
\int_{|x|>R}u_h+\int_{|x|>R}|u_{h+(g-h)_+}-u_f| \\ [4mm]
&\displaystyle\le \int_{|x|>R}u_h+\|(g-h)_+\|_{1} \le
\int_{|x|>R}u_h+\varepsilon.
\end{array}
$$
Hence, it is enough to consider the special case where $g$ is a regular,
radially symmetric, approximation of $M\chi_{\{|x|\le r_0\}}$, with support
contained in $\{|x|\le r_0+1\}$. Notice that the same property of radial
symmetry is true for the solution for any later time.

\noindent {\sc Step 3. }  We now write equation \eqref{eq-dep-cont} as
$$
u_n-\Delta z_n=g,\qquad z_n=(-\Delta)^{-(2-\sigma_n)/2}u_n^{m_n},
$$
whose weak formulation is
\begin{equation}
  \int_{\mathbb{R}^N} u_n\varphi+\int_{\mathbb{R}^N}\langle\nabla
  z_n,\nabla\varphi\rangle=\int_{\mathbb{R}^N} g\varphi
\label{weak-depcont}\end{equation} for every test function
$\varphi$. For radial solutions this equation reads (abusing notation)
$$
r^{N-1}u_n-(r^{N-1}z_n')'=r^{N-1}g\quad(=0\qquad \mbox{for } r> r_0+1).
$$

Since $u_n^{m_n}\in L^p({\mathbb{R}^N})$ for all $
\max\{1,1/m_n\}=\rho_n\le p\le\infty$, and since $N>2$ implies $m_n>(2-\sigma_n)/N$,
and thus $\rho_n\le N/(N-\sigma_n)$, applying \eqref{HLS} we get that
$z_n=(-\Delta)^{-(2-\sigma_n)/2}u_n^{m_n}\in L^r({\mathbb{R}^N})$ for every
$N/(N-2)\le r\le\infty$. On the other hand $z_n\in
C^{1,\gamma}({\mathbb{R}^N})$ for some $\gamma\in(0,1)$ depending on
$\widetilde\sigma=\max_n\{\sigma_n\}<2$
\cite{silvestre}. We even have, as in \cite{Cabre-Sire}, the
estimate $\|z_n\|_{C^{1,\gamma}}\le c$ independent on $n$. On the
other hand, the quantity $\Phi_n=r^{N-1}z_n'$ will be seen as a
flux, and indeed
 $$
\int_R^\infty r^{N-1}u_n(r)\,dr=|\Phi_n(R)-\Phi_n(\infty)| \qquad
\mbox{for every } R\ge r_0+1,
$$
which implies that the uniform control of the tails is equivalent to the
uniform control of $\Phi_n(r)$ for large $r$. It is clear that there exists a
limit flux, $\Phi_n(\infty)=\lim_{r\to\infty}\Phi_n(r)$, since it is monotone
increasing, $\Phi_n'(r)=r^{N-1}u_n\ge0$ for $r\ge r_0+1$. Notice also that
integrating the equation, we get
$$
\Phi_n(\infty)-\Phi_n(0)=\int_0^\infty
r^{N-1}u_n(r)\,dr-\int_0^\infty r^{N-1}g(r)\,dr\le0.
$$
Since $\Phi_n(0)=0$, we get $\Phi_n(\infty)\le0$. Moreover, $\Phi_n(r)\le0$ for
$r\ge r_0+1$.

Let us prove that $\Phi_n(\infty)=0$. If there exist constants
$C,r_1>0$ such that $\Phi_n(r)\le - C$ for $r\ge r_1$,   after an
integration we get
$$
z_n(r)\ge Cr^{2-N}\quad \mbox{for } r\ge r_0+1.
$$
This is a contradiction with the property $z_n\in L^{N/(N-2)}({\mathbb{R}^N})$.
Thus $\Phi_n(\infty)=0$.

\noindent {\sc Step 4. } We now have to estimate this limit more carefully. We
claim that $\Phi_n(r)$ is small for large $r$ uniformly in $n$,
which will give the desired uniform control of the tails.

Hardy-Littlewood-Sobolev's inequality~\eqref{HLS} gives a control of the
norm $\|z_n\|_{\frac N{N-2}}$ in terms of the norm $\|u_n^{m_n}\|_{\rho_n}$.
But this norm  can be estimated easily in terms of the norms $\|g\|_1$ and
$\|g\|_\infty$ and the value $\widetilde{m}=\max_n\{m_n\}$, with constants
independent on $n$. Therefore
\begin{equation}
\label{eq:estimate.gradients.z_n} \int_0^\infty r^{N-1}z_n^{\frac
N{N-2}}(r)\,dr\le C
\end{equation}
for every $n\ge1$.

We now adapt some ideas from \cite{Benilan-Crandall}. For
$\varepsilon>0$ given, if we take
$R(\varepsilon)=r_0e^{2C\varepsilon^{-\frac N{N-2}}}$, we have the
estimate
$$
\int_{r_0}^{R(\varepsilon)} r^{N-1}(\varepsilon r^{2-N})^{\frac N{N-2}}\,dr\ge 2C.
$$

Integrating the inequality
$$
z_n'(s)\le\left(\frac ts\right)^{N-1}z_n'(t), \qquad r_0+1\le r\le s,
$$
in $s$ for $s\in[r,t]$, and taking $t\ge2r$, we get, since
$z_n\ge0$,
$$
|\Phi_n(t)|=\left|t^{N-1}z_n'(t)\right|\le cr^{N-2}z_n(r);
$$
i.e., $
|\Phi_n(2R(\varepsilon))|r^{2-N}\le cz_n(r)$
whenever $R(\varepsilon)\ge r$. Putting all together we get
$$
\begin{array}{rl}
\displaystyle c\int_{r_0}^{R(\varepsilon)}
r^{N-1}(|\Phi_n(2R(\varepsilon))| r^{2-N})^{\frac N{N-2}}\,dr
&\displaystyle\le \int_{r_0}^{R(\varepsilon)} r^{N-1}z_n^{\frac
N{N-2}}(r)\,dr \\ [3mm] &\displaystyle\le C\le
\int_{r_0}^{R(\varepsilon)} r^{N-1}(\varepsilon r^{2-N})^{\frac
N{N-2}}\,dr,
\end{array}
$$
which implies
$$
|\Phi_n(R)|\le|\Phi_n(2R(\varepsilon))|\le c\,\varepsilon,\qquad
\mbox{for every } R\ge 2R(\varepsilon),\; n\ge 1.
$$
This ends the uniform control of the tails.

\noindent {\sc Step 5. }  Summing up, we have obtained that along some
subsequence, which we also call $\{(m_n,\sigma_n)\}$, the following
convergence holds
$$
u_n=u_{m_n,\sigma_n}\rightarrow\;u_*\qquad \mbox{in } L^1({\mathbb{R}^N}),
$$
for some function $u_*$. What is left is the identification of the limit, that
is $u_*=u_{\overline m,\overline\sigma}$.

The convergence $u_n\rightarrow u_*$ in $L^1({\mathbb{R}^N})$
implies the convergence $u_n^{m_n}\rightarrow\;u_*^{\overline m}$ in
$L^r({\mathbb{R}^N})$, $r=\max\{1,1/\min_n\{m_n\}\,\}$. Therefore,
using  inequality \eqref{HLS}, $z_n\rightarrow
z=(-\Delta)^{-(2-\overline\sigma)/2}u_*^{\overline m}$ in $L^{\frac
N{N-2}}({\mathbb{R}^N})$.

We now take $z_n$ as test function in~\eqref{weak-depcont}, use that
$u_n\ge0$, and apply H\"{o}lder's inequality with exponents
$p=N/(N-2)$, $q=N/2$,  and obtain, thanks
to~\eqref{eq:estimate.gradients.z_n}, a uniform control of the
gradients of $z_n$,
$$
\int_{\mathbb{R}^N} |\nabla z_n|^2\le \|g\|_{\frac
N2}\,\|z_n\|_{\frac{N}{N-2}}\le C.
$$
Hence, $\nabla z_n\rightharpoonup\nabla z$ in $L^2({\mathbb{R}^N})$.
All this is enough to pass to the limit in \eqref{weak-depcont} and
show that the limit $u_*$ is indeed $u_{\overline
m,\overline\sigma}$.~\qed

\subsection{$\sigma\to2$}
We now study the upper limit $\sigma\to2$.

\begin{Theorem}\label{th:cont.sigma2}
The map $S$ is also continuous at $\sigma=2$.
\end{Theorem}

\noindent{\it Proof. } Without loss of generality we keep $m$ fixed
and let $\sigma_n\to2^-$. As before, it is enough to prove the
convergence of the semigroup for bounded functions. The proof uses
the extension technique introducing the vertical variable $y>0$ and
considering the extended problem in the upper-half space
\eqref{pe:local}. The convergence, in $L^1_{\rm
loc}(\mathbb{R}^{N+1}_+)$, of the sequence $w_n=E_{\sigma_n}(u_n^m)$
to some function $w_*$, as well as that of $\nabla w_n$ to $\nabla
w_*$ in in compact sets, works as before.

To identify the limit
of the trace $u_*=(\Tr(w_*))^{1/m}$, we follow an idea from
\cite{Cabre-Sire}. We take a factorized test function, $\xi(x)\eta(y)$, where
$\eta$ is a cut-off function, $\eta(y)=1$ for $y\le 1$, $\eta(y)=0$ for $y\ge
2$, and $\xi$ is a usual test function in the $x$ variables.  Using that the
measures $(2-\sigma_n) y^{1-\sigma_n}dy$ are probability measures on $(0,1)$
converging (in the weak-$*$ sense of measures) to a Dirac measure $\delta_0$,
applied to the sequence $\varphi_n(y)=\eta(y)\nabla_xu_n^{m}(x,y)$, we finally
arrive to
$$
\int_{\mathbb{R}^N} u_*\xi\,dx+\int_{\mathbb{R}^N} \langle\nabla u_*^{
m},\nabla\xi\rangle\,dx=\int_{\mathbb{R}^N} g\xi\,dx,
$$
since $\lim_{\sigma\to2^-}\frac{\mu_{\sigma}}{2-\sigma}=1$. This is
the weak formulation of the equation $ u_*-\Delta u_*^m=g$.~\qed

\medskip

\section{Comments and extensions}\label{sec.com}
\setcounter{equation}{0}

\noindent {\sc $\bullet$ Limit.} The limit case $\sigma\to0$ looks very interesting. Notice that when $\sigma=0$, the critical exponent is $m_*=1$. Formally, the limit equation for $\sigma=0$ and $m\ge1$ is the ODE
\begin{equation}
\frac{\partial u}{\partial t}+ |u|^{m-1}u=0.
\label{ode}\end{equation}
In case the initial datum $u(x,0)=f(x)$ is a function defined pointwise, the
ODE can be explicitly solved, giving the formulas
\begin{equation}
u(x,t)=\left\{\begin{array}{ll}
(|f(x)|^{1-m}+(m-1)t)^{-1/(m-1)}\mbox{sign}(f(x)) &\quad \mbox{if} \ m>1, \\
f(x)\,e^{-t} &\quad \mbox{if} \ m=1.
\end{array}\right.
\label{sol-ode}\end{equation}
We see from these formulas that mass is not conserved; it  decays instead. Thus
we can not have convergence in $L^1(\mathbb{R}^N)$ as $\sigma\to0^+$. Observe
also by passing that when $m>1$ the decay rate $t^{-1/(m-1)}$ of
$\|u(\cdot,t)\|_\infty$ agrees with that of formula
\eqref{eq:L-inf-pp}.

On the other hand, the question of continuous
dependence below the critical exponent $m_*$, in some weighted norm,
is interesting and will be the subject of a future work.

 \noindent{\sc $\bullet$  Bounded domains. } We
have presented the basic facts for a theory in Sections
\ref{sec.bounded.domain} and  \ref{ssec.bounded}. We point out that
there are other ways to understand the Cauchy-Dirichlet problem in a
bounded domain with homogeneous \lq boundary data'. For example, one
may look for solutions to the problem
$$
\left\{\begin{array}{ll} \displaystyle\frac{\partial u}{\partial
t}=C_{N,\sigma }\mbox{
P.V.}\int_{\mathbb{R}^N}\frac{u^m(y,t)-u^m(x,t)}{|x-y|^{N+\sigma}}\,
dy, &\quad x\in \Omega,\, t>0 \\ [3mm] u(x,t)=0, &\quad
x\in\mathbb{R}^N\setminus\Omega,\, t\ge0,\\ [3mm] u(x,0)=f(x),
&\quad x\in \Omega,
\end{array}\right.
$$
which is different from Problem~\eqref{eq:main.bounded}. This approach has been recently used by Kim and Lee~\cite{Kim-Lee-1},~\cite{Kim-Lee-2}, who have addressed some
important issues in that framework, such as existence, regularity and
asymptotic behaviour.

The continuous dependence proofs of the last section extend to the case of solutions of Problem~\eqref{eq:main.bounded}.

\noindent{\sc $\bullet$ Related work. } Recently, Cifani and Jakobsen have
studied the existence of solutions of the diffusion-convection
equation $\partial_tu+\nabla \cdot f(u)+ (-\Delta)^{\sigma/2} A(u)
=0$, in the framework of (Kruzhkov-style) entropy solutions
\cite{CJ}. However, their assumptions on the nonlinearities
exclude~\eqref{eq:main}, unless $m=1$.

\vskip .5cm

{\small
\noindent \textsc{Acknowledgments.} All the authors supported by
Spanish Projects MTM2008-06326-C02-01 and -02. The last author partially supported
by MSRI, Berkeley, USA.   
We want to thank M. Bonforte and G. Grillo for some remarks on the
subcritical case.


}


\newpage

\noindent{\bf Addresses:}

\noindent{\sc A. de Pablo: } Departamento de Matem\'{a}ticas, Universidad Carlos III de Madrid, 28911 Legan\'{e}s,
Spain. (e-mail: arturo.depablo@uc3m.es).

\noindent{\sc F. Quir\'{o}s: }
Departamento de Matem\'{a}ticas, Universidad Aut\'{o}noma de Madrid, 28049 Madrid, Spain.
(e-mail: fernando.quiros@uam.es).

\noindent{\sc A. Rodr\'{\i}guez: }
Departamento de Matem\'{a}tica, ETS Arquitectura, Universidad Polit\'{e}cnica de Madrid, 28040 Madrid, Spain. (e-mail: ana.rodriguez@upm.es).

\noindent{\sc J.~L. V\'{a}zquez: }
Departamento de Matem\'{a}ticas, Universidad Aut\'{o}noma de Madrid, 28049
Madrid, Spain. (e-mail: juanluis.vazquez@uam.es). 

\end{document}